\documentclass[12pt]{article}

\title{Reverse Mathematics of Compact Countable Second-Countable Spaces}
\author{Fran{\c c}ois G. Dorais}
\date{\today}

\usepackage{verbatim}
\usepackage{enumerate}
\usepackage{amsmath}
\usepackage{amssymb}
\usepackage{amsthm}
\usepackage{hyperref}

\newcounter{count}[section]

\theoremstyle{plain}

\newtheorem{lemma}[count]{Lemma}
\newtheorem{proposition}[count]{Proposition}

\theoremstyle{definition}
\newtheorem{definition}[count]{Definition}
\newtheorem{example}[count]{Example}

\newtheorem{question}[count]{Question}
\theoremstyle{plain}
\newtheorem{claim}{Claim}[count]

\newcommand{\set}[1]{\lbrace#1\rbrace}
\newcommand{\seq}[1]{\langle#1\rangle}

\newcommand{\lthen}{\rightarrow}
\newcommand{\liff}{\leftrightarrow}

\newcommand{\cat}{{{}^\frown}}
\newcommand{\res}{{\upharpoonright}}
\newcommand{\rem}{\setminus}

\newcommand{\N}{\mathbb{N}}
\newcommand{\Z}{\mathbb{Z}}

\newcommand{\RCA}{\ensuremath{\mathsf{RCA}_0}}
\newcommand{\WKL}{\ensuremath{\mathsf{WKL}_0}}
\newcommand{\ACA}{\ensuremath{\mathsf{ACA}_0}}
\newcommand{\Ind}[1]{\ensuremath{\mathsf{I}{#1}}}
\newcommand{\Bnd}[1]{\ensuremath{\mathsf{B}{#1}}}

\begin{document}

\maketitle

\begin{abstract}\noindent
  We study the reverse mathematics of the theory of countable second-countable topological spaces, with a focus on compactness.
  We show that the general theory of such spaces works as expected in the subsystem \ACA\ of second-order arithmetic, but we find that many unexpected pathologies can occur in weaker subsystems.
  In particular, we show that \RCA\ does not prove that every compact discrete countable second-countable space is finite and that \RCA\ does not prove that the product of two compact countable second-countable spaces is compact.
  To circumvent these pathologies, we introduce strengthened forms of compactness, discreteness, and Hausdorffness which are better behaved in subsystems of second-order arithmetic weaker than \ACA.
\end{abstract}

\section{Introduction}

The goal of reverse mathematics is to understand the various set comprehension principles that are necessary to prove central theorems of everyday mathematics.
The context for this analysis is usually second-order arithmetic.
This approach has been extraordinarily successful for the analysis of major results in algebra, analysis, and combinatorics.
General topology presents a unique challenge for this approach since the very definition of topological spaces involves three layers of types: points, sets of points, and families of sets of points.
Thus, there is no hope for second-order arithmetic to fully grasp the rich variety of topological spaces.

Nevertheless, considerable efforts have been made to understand a wide spectrum of topological spaces in the context of second-order arithmetic.
Notably, Mummert~\cite{Mummert06} studied the class of MF spaces.
This class of spaces is broad enough to include all complete metric spaces as well as many non-metrizable spaces.
However, MF spaces are still limited in many ways, for example all such spaces are $T_1$, second-countable, and strong Choquet~\cite{MummertStephan10}.
Another approach was used by Hunter~\cite{Hunter08} who studied topological spaces in the context of higher-order arithmetic.
This approach is interesting since the presence of higher types allows for a direct understanding of the topology on a space.
However, this breaks with the longstanding tradition of using second-order arithmetic for reverse mathematics and thus makes difficult the comparison of theorems in general topology and those in other fields of mathematics.

In this paper, we will focus on the reverse mathematics of countable second-countable spaces.
These have the advantage that they can be directly encoded in the classical context of second-order arithmetic.
Indeed, the set of points can be identified with a subset of the natural numbers, and a base for such a space can be explicitly listed as an indexed sequence of sets of natural numbers.
Although the full topology is still not directly accessible, common topological concepts can be understood as properties of basic open sets.
We will investigate several such concepts, e.g., compactness, discreteness, Hausdorffness.
We will show that all of these concepts work as expected in the susbsystem \ACA, but that many unexpected phenomena can occur in weaker subsystems.
For example:
\begin{itemize}
\item \RCA\ does not prove that every compact discrete countable second-countable space is finite.
  (Example~\ref{E:DiscreteCompact}.)
\item \RCA\ does not prove that every compact countable second-countable space is sequentially compact.
  (Example~\ref{E:DiscreteCompact}.)
\item \RCA\ does not prove that the order topology of a countable complete linear ordering is compact.
  (Example~\ref{E:DiscreteComplete}.)
\item \RCA\ does not prove that the product of two compact countable second-countable spaces is compact.
  (Example~\ref{E:Tychonoff}.)
\end{itemize}
However, in each case, we identify one or more strengthened forms of the concept that can be used in \RCA\ with fewer unexpected consequences.

A standard reference for subsystems of second-order arithmetic and their use in reverse mathematics is Simpson~\cite{Simpson09}.
Formal definitions and properties of the basic subsystems \RCA, \WKL, and \ACA\ can be found there.
Due to the highly combinatorial nature of our subject matter, we will make use of several conventions throughout this paper.
The goal of these conventions is to minimize the amount of arithmetical coding in the statement and proofs of our results.

\paragraph{Finite sequences.}
Finite sequences are coded using natural numbers in such a way that all the basic operations on finite sequences are primitive recursive.
There are many ways to accomplish this; see Smory{\'n}ski~\cite{Smorynski77} for a detailed account.
Since the details of this coding are immaterial, we will not impose one particular choice on the reader.

We will use the notation \(X^{<\infty}\) to denote the set of all sequences of elements of the set \(X.\)
The length of a sequence \(x\) will be denoted \(|x|.\)
For a natural number \(n,\) we will write \(X^n\) for the subset of \(X^{<\infty}\) consisting of sequences of length \(n.\)
We will write \(\seq{x_0,\dots,x_{\ell-1}}\) for the sequence of length \(\ell\) whose \(i\)-th element is \(x_{i-1}.\)

\paragraph{Finite sets.}
Finite sets are identified with the finite sequences which enumerate them in increasing order.
We will write \(X^{[<\infty]}\) for the set of all finite sets of elements of the set \(X.\)
We will write \(\set{x_0,\dots,x_{\ell-1}}\) for the set whose elements are enumerated by the sequence \(\seq{x_0,\dots,x_{\ell-1}}\) (not necessarily in increasing order).

We will handle finite sets as if they were plain sets.
For example, we will write \(x \in a\) to abbreviate \((\exists i < |a|)(x = a_i)\) and we will write \(a \subseteq b\) to abbreviate \((\forall i < |a|)(\exists j < |b|)(a_i = b_j).\)
We will freely take unions and intersections of finite sets since these are primitive recursive operations.

\paragraph{Enumerable sets.}
An \emph{enumerable set} is a nondecreasing sequence \(A = \seq{a_n}_{n=0}^\infty\) of finite sets.
This sequence is intended to represent the union \(\bigcup_{n=0}^\infty a_n,\) which does not provably exist in \RCA.

We will handle enumerable sets as if they were plain sets.
For example, we will write \(x \in A\) to abbreviate \((\exists n)(x \in a_n)\) and we will write \(A \subseteq B\) to abbreviate \((\forall x)(x \in A \lthen x \in B).\) 
Although this is somewhat confusing, we will also write \(A = B\) to abbreviate \((\forall x)(x \in A \liff x \in B)\) since the pointwise equality of the sequences \(A\) and \(B\) will never be relevant.
We will freely take countable unions and finite intersections of enumerable sets since these can be described in a canonical way.

Note that for every \(\Sigma^0_1\) formula \(\phi(x)\) there is an enumerable set \(A\) such that \((\forall x)(x \in A \liff \phi(x));\) we will abbreviate this by writing \(A = \set{x : \phi(x)}.\)
In fact, we will often use set comprehension to define enumerable sets.
In other words, we will often write \(A = \set{x : \phi(x)}\) to abbreviate the routine process by which \(\phi(x)\) is used to construct a nondecreasing sequence \(A = \seq{a_n}_{n=0}^\infty\) of finite sets such that \((\exists n)(x \in a_n) \liff \phi(x).\)

\section{Countable Second-Countable Spaces}

We begin our study of countable second-countable spaces by laying out some definitions that will be used throughout this paper.

\begin{definition}[\RCA]
  A \emph{base} for a topology on a set \(X\) is an indexed sequence \(\mathcal{U} = \seq{U_i}_{i \in I}\) of subsets of \(X\) together with a function \(k: X \times I \times I \to I\) such that the following two properties hold.
  \begin{itemize}
  \item If \(x \in X\) then \(x \in U_i\) for some \(i \in I.\)
  \item If \(x \in U_i \cap U_j\) then \(x \in U_{k(x,i,j)} \subseteq U_i \cap U_j.\)
  \end{itemize}
\end{definition}

\noindent
Over \ACA, we can omit the function \(k\) and relax the second property to:
\begin{itemize}
\item If \(x \in U_i \cap U_j\) then \(x \in U_k \subseteq U_i \cap U_j\) for some \(k \in I.\)
\end{itemize}
Indeed, a suitable function always exists since the property required of \(k\) is arithmetical.

\begin{definition}[\RCA]
  A \emph{countable second-countable space} is a triple \((X,\mathcal{U},k)\) where \(\mathcal{U} = \seq{U_i}_{i \in I}\) and \(k: X \times I \times I \to I\) form a base for a topology on the set \(X.\) 
\end{definition}

\noindent
Over \ACA, we usually omit the third parameter \(k\) for the reasons explained above.

When defining open sets in a countable second-countable space \((X,\mathcal{U},k),\) we bump into the difficulty that not all internal unions of basic open sets need to exist in \RCA.
For this reason, we resort to codes for open sets.

\begin{definition}[\RCA]
  Let \((X,\mathcal{U},k)\) be a countable second-countable space where \(\mathcal{U} = \seq{U_i}_{i \in I}.\)
  An \emph{open code} in \((X,\mathcal{U},k)\) is an enumerable subset \(A\) of~\(I.\)
\end{definition}

\noindent
The open code \(A\) is intended to represent the enumerable set \(\bigcup_{i \in A} U_i,\) which may not be a set in a model of \RCA.

\begin{definition}[\RCA]
  Let \((X,\mathcal{U},k)\) be a countable second-countable space.
  A subset \(G\) of \(X\) is \emph{effectively open} if there is an open code \(A\) such that \(G = \bigcup_{i \in A} U_i.\)
  Dually, a subset \(F\) of \(X\) is \emph{effectively closed} if the complement of \(F\) in \(X\) is effectively open.  
\end{definition}

\noindent
Note that this is indeed an effective version of the local definition of openness.
Namely, if \(G\) is effectively open and \(A\) is an open code such that \(G = \bigcup_{i \in A} U_i,\) then for every \(x \in G\) we can effectively search through \(A\) for an index \(i\) such that \(x \in U_i.\)
Therefore, there is a function \(g:G \to A\) such that \(x \in U_{g(x)} \subseteq G\) for every \(x \in G.\)

Let \((X,\mathcal{U},k)\) and \((Y,\mathcal{V},\ell)\) be countable second-countable spaces.
The continuity of a function \(f:X \to Y\) is defined locally.
The function \(f\) is \emph{continuous at the point \(x \in X\)} if for every basic open neighborhood \(V_j\) of \(f(x)\) there is a basic open neighborhood \(U_i\) of \(x\) such that \(U_i \subseteq f^{-1}[V_j];\) the function \(f\) is \emph{continuous} if it is continuous at every point of \(X.\)
This definition of continuity makes perfect sense in \ACA, but it is impractical in \RCA\ since there may be no way to effectively compute an index for \(U_i\) given \(x\) and an index for \(V_j.\)

\begin{definition}[\RCA]
  Let \((X,\mathcal{U},k)\) and \((Y,\mathcal{V},\ell)\) be countable second-countable spaces where \(\mathcal{U} = \seq{U_i}_{i \in I}\) and \(\mathcal{V} = \seq{V_j}_{j \in J}.\)
  A function \(f:X \to Y\) is \emph{effectively continuous} if there is a function \(\phi:X \times J \to I\) such that if \(f(x) \in V_j\) then \(x \in U_{\phi(x,j)} \subseteq f^{-1}[V_j].\)
\end{definition}

\noindent
With this definition of effective continuity, we can show that preimages of open sets are open.

\begin{proposition}[\RCA]\label{P:ContinuityRCA}
  Let \((X,\mathcal{U},k)\) and \((Y,\mathcal{V},\ell)\) be countable second-countable spaces and let \(f:X \to Y\) be an effectively continuous function.
  If \(B\) is an open code in \((Y,\mathcal{V},\ell)\) then we can find an open code \(A\) in \((X,\mathcal{U},k)\) such that \(\bigcup_{i \in A} U_i = f^{-1}[\bigcup_{j \in B} V_j].\)
\end{proposition}

\begin{proof}
  Write \(\mathcal{U} = \seq{U_i}_{i \in I},\) \(\mathcal{V} = \seq{V_j}_{j \in J},\) and let \(\phi:X \times J \to I\) witness the effective continuity of \(f.\)
  Given an enumerable set \(B \subseteq J,\) the enumerable set \[A = \set{\phi(x,j) : j \in B \land f(x) \in V_j}\] is as required.
\end{proof}

\begin{proposition}[\RCA]
  Let \((X,\mathcal{U},k)\) and \((Y,\mathcal{V},\ell)\) be countable second-countable spaces and let \(f:X \to Y\) be an effectively continuous function.
  If \(V\) is an effectively open subset of \(Y,\) then \(f^{-1}[V]\) is an effectively open subset of \(X.\)
\end{proposition}

\noindent
Thus, over \ACA, we recover the usual global definition of continuity: a function \(f:X \to Y\) is continuous if and only if \(f^{-1}[V]\) is open in \((X,\mathcal{U})\) for every open set \(V\) in \((Y,\mathcal{V}).\)

As usual, an \emph{\textup(effective\textup) homeomorphism} between two spaces \((X,\mathcal{U},k)\) and \((Y,\mathcal{V},\ell)\) is a bijection \(f:X \to Y\) such that both \(f\) and \(f^{-1}\) are \textup(effectively\textup) continuous.

\begin{definition}[\RCA]
  Let \(\mathcal{U},k\) and \(\mathcal{V},\ell\) be two countable bases on the set \(X.\)
  We say that \(\mathcal{U},k\) and \(\mathcal{V},\ell\) are \emph{effectively equivalent} if the identity function on \(X\) is an effective homeomorphism between \((X,\mathcal{U},k)\) and \((X,\mathcal{V},\ell).\)
\end{definition}

\noindent
This definition will be useful to identify when a property of a space is a topological property and not a property of the presentation of a space: truly topological properties should be invariant under effectively equivalent bases for the same topology.

Subspaces and product spaces are defined as expected.

\begin{definition}[\RCA]
  Suppose \((X,\mathcal{U},k)\) is a countable second-countable space with \(\mathcal{U} = \seq{U_i}_{i \in I}.\) 
  If \(X' \subseteq X\) then the corresponding \emph{subspace} \((X',\mathcal{U}',k')\) is defined by \(U'_i = U_i \cap X'\) for all \(i \in I,\) and \(k' = k\res(X'\times I \times I).\)
\end{definition}

\begin{definition}[\RCA]
  Suppose \((X,\mathcal{U},k)\) and \((Y,\mathcal{V},\ell)\) are countable second-countable space with \(\mathcal{U} = \seq{U_i}_{i \in I}\) and \(\mathcal{V} = \seq{V_j}_{j \in J}.\)
  The \emph{product space} \((X,\mathcal{U},k)\times(Y,\mathcal{V},\ell)\) is the countable second-countable space \((X \times Y, \mathcal{W}, m),\) where \(\mathcal{W} = \seq{U_i \times V_j}_{\seq{i,j} \in I \times J}\) and \[m(\seq{x,y},\seq{i,j},\seq{i',j'}) = \seq{k(x,i,i'),\ell(y,j,j')}.\]
\end{definition}

\noindent
It is straightforward to check that these two constructs obey the definition of countable second-countable space from above.

A common way to construct a topological space is to specify a subbase for the topology.
We will use this very often in this paper, so we introduce some formal notation for going from a subbase to a base for a countable second-countable topological space.

\begin{definition}[\RCA]
  If \(\mathcal{B} = \seq{B_i}_{i \in I}\) is a sequence of subsets of \(X,\) then we define \(\mathcal{B}^* = \seq{B^*_s}_{s \in I^{[<\infty]}}\) by \(B^*_s = \bigcap_{i \in s} B_i,\) with the convention that \(B^*_\varnothing = X.\)
\end{definition}

\noindent
Note that the sequence \(\mathcal{B}^*\) is easily definable by primitive recursion on notation, so there is no trouble forming \(\mathcal{B}^*\) in models of \RCA.

\begin{proposition}[\RCA]\label{P:Subbase}
  If \(\mathcal{B} = \seq{B_i}_{i \in I}\) is an arbitrary sequence of subsets of \(X,\) then \((X,\mathcal{B}^*,k^*)\) is a countable second-countable space, where \(k^*:X \times I^{[<\infty]} \times I^{[<\infty]} \to I^{[<\infty]}\) is invariably defined by \(k^*(x,s,t) = s \cup t.\)
\end{proposition}

\noindent
Many of the examples of countable second-countable spaces that follow will be defined by specifying a subbase \(\mathcal{B}\) for the open sets and defining the base as \(\mathcal{B}^*.\)

The following definition will be useful in our study of compactness.

\begin{definition}[\RCA]
  A sequence \(\seq{U_i}_{i \in I}\) of subsets of a set \(X\) has a \emph{finite cover relation} if there exists a set \(C \subseteq I^{[<\infty]}\) such that \[\set{i_0,\dots,i_{\ell-1}} \in C \liff X = U_{i_0} \cup\cdots\cup U_{i_{\ell-1}}.\]
\end{definition}

\noindent
Of course, over \ACA, every sequence of sets has a finite cover relation.

The following result will be helpful in the construction of countable second-countable spaces whose bases have finite cover relations.

\begin{lemma}[\RCA]\label{L:SubFiniteCover}
  If \(\mathcal{B}\) is a sequence of subsets of \(X\) with a finite cover relation, then \(\mathcal{B}^*\) has a finite cover relation too.
\end{lemma}

\begin{proof}
  Write \(\mathcal{B} = \seq{B_i}_{i \in I}\) and let \(C \subseteq I^{[<\infty]}\) be the finite cover relation for \(\mathcal{B}.\)
  Define \(C^* \subseteq (I^{[<\infty]})^{[<\infty]}\) by \(\set{s_0,\dots,s_{\ell-1}} \in C^*\) if and only if \(\set{i_0,\dots,i_{\ell-1}} \in C\) for every \(\seq{i_0,\dots,i_{\ell-1}} \in s_0 \times\cdots\times s_{\ell-1}.\)
  This is correct because \[\bigcup_{n=0}^{\ell-1} B^*_{s_n} = \bigcup_{n=0}^{\ell-1} \bigcap_{i \in s_n} B_i = \bigcap_{\seq{i_0,\dots,i_{\ell-1}} \in s_0\times\cdots\times s_{\ell-1}} \bigcup_{n=0}^{\ell-1} B_{i_n},\] with the usual convention that empty intersections equal \(X.\)
\end{proof}

\section{Compactness}

For compactness of countable second-countable spaces, we would like to use the traditional definition: every open cover has a finite subcover.
However, arbitrary open covers are third-order objects which are inaccessible to second-order arithmetic.
Therefore, we use the following variant which is equivalent to the usual definition and avoids third-order objects.

\begin{definition}[\RCA]
  Let \((X,\mathcal{U},k)\) be a countable second-countable space with \(\mathcal{U} = \seq{U_i}_{i \in I}.\)
  We say that \((X,\mathcal{U},k)\) is \emph{compact} if for every enumerable set \(A \subseteq I\) such that \(X = \bigcup_{i \in A} U_i\) there is a finite set \(a \subseteq A\) such that \(X = \bigcup_{i \in a} U_i.\)
\end{definition}

\noindent
A priori, compactness of a countable second-countable space depends on the choice of base for the topology.
However, the next proposition shows that this is not the case.

\begin{proposition}[\RCA]\label{P:CompactImage}
  Let \((X,\mathcal{U},k)\) and \((Y,\mathcal{V},\ell)\) be countable second-countable spaces.
  If there is an effectively continuous surjection \(f:X \to Y\) and \((X,\mathcal{U},k)\) is compact then \((Y,\mathcal{V},\ell)\) is compact too.
\end{proposition}

\begin{proof}
  Write \(\mathcal{U} = \seq{U_i}_{i \in I}\) and \(\mathcal{V} = \seq{V_j}_{j \in J},\) and let \(\phi:X \times J \to I\) witness the effective continuity of the surjection \(f:X \to Y.\)
  Given an enumerable set \(B \subseteq J\) such that \(Y = \bigcup_{j \in B} V_j,\) consider the enumerable set \[A = \set{\phi(x,j) : j \in B \land f(x) \in V_j}.\]
  Note that \(X = \bigcup_{i \in A} U_i.\)
  Since \((X,\mathcal{U},k)\) is compact, there is a finite set \(a \subseteq A\) such that \(X = \bigcup_{i \in a} U_i.\)
  By definition of \(A,\) there is a finite set \(b \subseteq B\) such that \[a \subseteq \set{\phi(x,j) : j \in b \land f(x) \in V_j}.\]
  We claim that \(Y = \bigcup_{j \in b} V_j.\)
  Given \(y_0 \in Y,\) first find \(x_0 \in X\) such that \(f(x_0) = y_0\) and then find \(i \in a\) such that \(x_0 \in U_i.\)
  Note that \(i = \phi(x',j)\) for some \(j \in b\) and some \(x \in X\) such that \(f(x) \in V_j.\)
  By definition of \(\phi,\) we then have \(U_i \subseteq f^{-1}[V_j]\) which means that \(y_0 = f(x_0) \in V_j.\)
\end{proof}

\noindent
It follows that if \(\mathcal{U},k\) and \(\mathcal{V},\ell\) are two effectively equivalent bases on the set \(X,\) then \((X,\mathcal{U},k)\) is compact if and only if \((Y,\mathcal{V},\ell)\) is compact.

While perfectly meaningful, the above definition of compactness is somewhat impractical in \RCA.
Indeed, there may be no effective way to determine which finite collections of basic open sets form finite covers.
Thus, a more effective version of compactness is the following.

\begin{definition}[\RCA]
  A countable second-countable space \((X,\mathcal{U},k)\) is \emph{effectively compact} if it is compact and its base \(\mathcal{U}\) has a finite cover relation.
\end{definition}

\noindent
Somewhat surprisingly, effective compactness also turns out to be a topological property.

\begin{proposition}[\RCA]\label{P:EffectivelyCompactImage}
  Let \((X,\mathcal{U},k)\) and \((Y,\mathcal{V},\ell)\) be countable second-countable spaces.
  If there is an effectively continuous surjection \(f:X \to Y\) and \((X,\mathcal{U},k)\) is effectively compact then \((Y,\mathcal{V},\ell)\) is effectively compact too.
\end{proposition}

\begin{proof}
  Write \(\mathcal{U} = \seq{U_i}_{i \in I}\) and \(\mathcal{V} = \seq{V_j}_{j \in J}.\) 
  Let \(C \subseteq I^{[<\infty]}\) be the finite cover relation for \(\mathcal{U}\) and let \(\phi: X \times J \to I\) witness the effective continuity of the surjection \(f: X \to Y.\)

  The finite cover relation for \(\mathcal{V}\) is always defined by a \(\Pi^0_1\) statement.
  Using \(C,\) \(f,\) and \(\phi,\) we can also describe it by a \(\Sigma^0_1\) statement:
  \[\bigcup_{j \in t} V_t = Y \liff (\exists s \in C)(\forall i \in s)(\exists x \in X)(\exists j \in t)(f(x) \in V_j \land i = \phi(x,j)).\]
  Therefore, \(\Delta^0_1\) comprehension suffices to comprehend the finite cover relation for \(\mathcal{V}.\)
\end{proof}

\noindent
It follows that if \(\mathcal{U},k\) and \(\mathcal{V},\ell\) are two effectively equivalent bases on the set \(X,\) then \((X,\mathcal{U},k)\) is effectively compact if and only if \((Y,\mathcal{V},\ell)\) is effectively compact.

The next example shows that \ACA\ is the minimal subsystem of second-order arithmetic in which all countable second-countable compact spaces are effectively compact.

\begin{example}\label{E:Compact}
  \emph{A compact space which is not effectively compact from a failure of \(\Pi^0_1\) comprehension.}

  Let \(\phi(i)\) be any \(\Pi^0_1\) formula and write \(\phi(i) \liff (\forall x)\phi_0(i,x)\) where \(\phi_0(i,x)\) is bounded.
  Consider the countable second-countable space \((\N,\mathcal{B}^*),\) where the subbase \(\mathcal{B} = \seq{B_i}_{i \in \N}\) consists of the sets \[B_i = \set{x \in \N : \phi_0(i,x) \lor (\exists y < x)\lnot\phi_0(i,y)}.\]
  Note that \((\N,\mathcal{B}^*)\) is a compact space since all nonempty basic open sets are cofinite.
  Suppose that \(C\) is a finite cover relation for \(\mathcal{B}^*,\) then \(\phi(i) \liff \set{\set{i}} \in C.\)
  Therefore, if \(\phi(i)\) witnesses the failure of \(\Pi^0_1\) comprehension, then \((\N,\mathcal{B}^*)\) is compact but not effectively compact.
\end{example}

Effectively closed subspaces of compact and effectively compact spaces behave as expected in \RCA.

\begin{proposition}[\RCA]\label{P:ClosedCompact}
  Let \((X,\mathcal{U},k)\) be a countable second-countable compact space.
  Every effectively closed subspace of \((X,\mathcal{U},k)\) is compact.
\end{proposition}

\begin{proof}
  Suppose \(X_0 \subseteq X\) is effectively closed.
  Say, \(A_1 \subseteq I\) is an enumerable set such that \(X \rem X_0 = \bigcup_{i \in A_1} U_i.\)
  
  Given an enumerable set \(A_0 \subseteq I\) such that \(X_0 \subseteq \bigcup_{i \in A_0} U_i,\) we have \(X = \bigcup_{i \in A_0 \cup A_1} U_i.\) 
  Since \((X,\mathcal{U},k)\) is compact, there is a finite set \(a \subseteq A_0 \cup A_1\) such that \(X = \bigcup_{i \in a} U_i.\)
  Since \(X_0 \cap U_i = \varnothing\) for every \(i \in A_1,\) we necessarily have \(X_0 \subseteq \bigcup_{i \in a_0} U_i\) where \(a_0 = a \cap A_0.\)
\end{proof}

\begin{proposition}[\RCA]\label{P:ClosedEffCompact}
  Let \((X,\mathcal{U},k)\) be a countable second-countable effectively compact space.
  Every effectively closed subspace of \((X,\mathcal{U},k)\) is effectively compact.
\end{proposition}

\begin{proof}
  Let \(C\) be the finite cover relation for \((X,\mathcal{U},k).\)

  Suppose \(X_0 \subseteq X\) is effectively closed.
  Say, \(A_1 \subseteq I\) is an enumerable set such that \(X \rem X_0 = \bigcup_{i \in A_1} U_i.\)
  We know from Proposition~\ref{P:ClosedCompact} that \(X_0\) is relatively compact in \((X,\mathcal{U},k),\) so it suffices to show that \(X_0\) has a finite cover relation.

  The finite cover relation for \(X_0\) always admits a \(\Pi^0_1\) description.
  We can also give a \(\Sigma^0_1\) description as follows:
  \[X_0 \subseteq \bigcup_{i \in t} U_i \liff (\exists s \in A_1^{[<\infty]})(t \cup s \in C).\]
  Thus, \(\Delta^0_1\) comprehension suffices to comprehend the finite cover relation for \(X_0.\)
\end{proof}

\section{Sequential Compactness}

Let \((X,\mathcal{U},k)\) be a countable second-countable space.
A point \(x \in X\) is said to be an \emph{accumulation point} of the sequence \(\seq{x_n}_{n=0}^\infty\) of points of \(X\) if the set \[\set{n \in \N : x_n \in U_i}\] is infinite for every basic open neighborhood \(U_i\) of \(x.\)
We say that \((X,\mathcal{U},k)\) is \emph{sequentially compact} if every sequence of points of \(X\) has an accumulation point.

\begin{proposition}[\RCA]\label{P:SeqCompactCompact}
  Every countable second-countable sequentially compact space is compact.
\end{proposition}

\begin{proof}
  Let \(A = \seq{a_n}_{n=0}^\infty\) be an enumerable subset of \(I\) such that \(\bigcup_{i \in a_n} U_i \neq X\) for each \(n.\)
  Define \(x_n\) to be the least element of \(X \rem \bigcup_{i \in a_n} U_i.\)
  Let \(x\) be an accumulation point of \(\seq{x_n}_{n=0}^\infty.\)
  We claim that \(x \notin \bigcup_{i \in A} U_i.\)
  Indeed, suppose instead that \(x \in U_i\) for some \(i \in A.\)
  Choose \(n_0\) such that \(i \in a_{n_0},\) then \(x_n \notin U_i\) for all \(n \geq n_0,\) which is impossible since the set \(\set{n : x_n \in U_i}\) must be infinite.
\end{proof}

\noindent
Example~\ref{E:DiscreteCompact} shows that \ACA\ is necessary to prove the converse of Proposition~\ref{P:SeqCompactCompact}, even if compactness is strengthened to effective compactness.
However, the proposition does admit a partial converse if compactness is strengthened in a different way.

\begin{definition}[\RCA]
  A sequence \(\seq{U_i}_{i \in I}\) of subsets of \(X\) is \emph{\(\Sigma^0_n\)-compact} if for every \(\Sigma^0_n\) formula \(A(i)\) such that \[(\forall x \in X)(\exists i \in I)(A(i) \land x \in U_i)\] there is a finite set \(a \subseteq I\) such that \((\forall i \in a)A(i)\) and \(X = \bigcup_{i \in a} U_i.\)
  A countable second-countable space \((X,\mathcal{U},k)\) is \emph{\(\Sigma^0_n\)-compact} if the sequence \(\mathcal{U}\) is \(\Sigma^0_n\)-compact.
\end{definition}

\noindent
Note that \(\Sigma^0_1\)-compactness is equivalent to plain compactness.

The following proposition shows that \(\Sigma^0_n\)-compactness is remarkably robust when \(n \geq 2.\)
Indeed, one consequence is that \(\Sigma^0_n\)-compactness is invariant under all homeomorphisms, not just the effective ones.

\begin{proposition}[\RCA\ + \Ind{\Sigma^0_n}, \(n \geq 2\)]\label{P:NCompactImage}
  Let \((X,\mathcal{U},k)\) and \((Y,\mathcal{V},\ell)\) be countable second-countable spaces.
  If there is a continuous surjection \(f:X \to Y\) and \((X,\mathcal{U},k)\) is \(\Sigma^0_n\)-compact then \((Y,\mathcal{V},\ell)\) is \(\Sigma^0_n\)-compact too.
\end{proposition}

\begin{proof}
  Write \(\mathcal{U} = \seq{U_i}_{i \in I}\) and \(\mathcal{V} = \seq{V_j}_{j \in J}.\)
  Given a \(\Sigma^0_n\) formula \(B(j)\) such that \[(\forall y \in Y)(\exists j \in J)(B(j) \land y \in V_j),\] consider the formula \(A(i)\) defined by \[(\exists j \in J)(B(j) \land U_i \subseteq f^{-1}[V_j]).\]
  This is a \(\Sigma^0_n\) formula since \(n \geq 2,\) and note that \[(\forall x \in X)(\exists i \in I)(A(i) \land x \in U_i).\]
  Since \((X,\mathcal{U},k)\) is \(\Sigma^0_n\)-compact, there is a finite set \(a \subseteq I\) such that \((\forall i \in a)A(i)\) and \(X = \bigcup_{i \in a} U_i.\)
  By \Ind{\Sigma^0_n} and the definition of \(A(i),\) there is a finite set \(b \subseteq J\) such that \((\forall j \in b)B(j)\) and \[(\forall i \in a)(\exists j \in b)(U_i \subseteq f^{-1}[V_j]).\]
  We claim that \(Y = \bigcup_{j \in b} V_j.\)
  Given \(y_0 \in Y,\) first find \(x_0 \in X\) such that \(f(x_0) = y_0\) and then find \(i \in a\) such that \(x_0 \in U_i.\)
  Note that \(U_i \subseteq f^{-1}[V_j]\) for some \(j \in b,\) which means that \(y_0 = f(x_0) \in V_j.\)
\end{proof}

\noindent
Note that the continuous surjection is not required to be effectively continuous in Proposition~\ref{P:NCompactImage}.

Closed subspaces of \(\Sigma^0_n\)-compact spaces are also well behaved.
 
\begin{proposition}[\RCA\ + \Ind{\Sigma^0_n}, \(n \geq 2\)]\label{P:ClosedNCompact}
  Let \((X,\mathcal{U},k)\) be a countable second-countable \(\Sigma^0_n\)-compact space.
  Every closed subspace of \((X,\mathcal{U},k)\) is \(\Sigma^0_n\)-compact.
\end{proposition}

\begin{proof}
  Suppose \(X_0 \subseteq X\) is closed.  
  Given a \(\Sigma^0_n\) formula \(A_0(i)\) such that \[(\forall x \in X_0)(\exists i \in I)(A_0(i) \land x \in U_i),\] we have \[(\forall x \in X)(\exists i \in I)((A_0(i) \lor U_i \cap X_0 = \varnothing) \land x \in U_i).\] 
  The formula \(A_0(i) \lor U_i \cap X_0 = \varnothing\) is \(\Sigma^0_n\) formula since \(n \geq 2.\)
  Since \((X,\mathcal{U},k)\) is \(\Sigma^0_n\)-compact, there is a finite set \(a \subseteq I\) such that \((\forall i \in a)(A_0(i) \lor U_i \cap X_0 = \varnothing)\) and \(X = \bigcup_{i \in a} U_i.\)
  We then necessarily have \(X_0 \subseteq \bigcup_{i \in a_0} U_i\) where \(a_0 = \set{i \in a : A_0(i)},\) which exists by \Ind{\Sigma^0_n}.
\end{proof}

\noindent
Note again that the closed subspace is not required to be effectively closed in Proposition~\ref{P:ClosedNCompact}.

The next proposition is useful for producing examples of \(\Sigma^0_n\)-compact spaces. 
Indeed, it shows that the non-effectively compact space from Example~\ref{E:Compact} is \(\Sigma^0_n\)-compact under \Ind{\Sigma^0_n}.
Thus \(\Sigma^0_n\)-compactness is a strengthening of compactness which is orthogonal to effective compactness.
In fact, it makes sense to talk about effectively \(\Sigma^0_n\)-compact spaces.

\begin{proposition}[\RCA\ + \Ind{\Sigma^0_n}]\label{P:CofiniteNCompact}
  Let \((X,\mathcal{U},k)\) be a countable second-countable space wherein some point \(x_0\) is such that all basic open neighborhoods of \(x_0\) are cofinite, then \((X,\mathcal{U},k)\) is \(\Sigma^0_n\)-compact.
\end{proposition}

\begin{proof}
  Suppose that \(A(i)\) is a \(\Sigma^0_n\) formula such that \[(\forall x \in X)(\exists i \in I)(A(i) \land x \in U_i).\]
  A fortiori, there is an \(i_0 \in I\) such that \(A(i_0)\) and \(x_0 \in U_{i_0}.\)
  Since \(U_{i_0}\) is cofinite, let \(\seq{x_1,\dots,x_{\ell-1}}\) enumerate the complement of \(U_{i_0}.\)
  By \Ind{\Sigma^0_n}, there is a corresponding sequence \(\seq{i_1,\dots,i_{\ell-1}}\) of elements of \(I\) such that \(A(i_m)\) and \(x_m \in U_{i_m}\) for each \(m < \ell.\)
  The set \(a = \set{i_0,i_1,\dots,i_{\ell-1}}\) is then as required for \(\Sigma^0_n\)-compactness.
\end{proof}

We are now ready to prove the promised partial converse of Proposition~\ref{P:SeqCompactCompact}.

\begin{proposition}[\RCA\ + \Bnd{\Sigma^0_2}]\label{P:2CompactSeqCompact}
  Every countable second-countable \(\Sigma^0_2\)-compact space is sequentially compact.
\end{proposition}

\begin{proof}
  Let \((X,\mathcal{U},k)\) be a countable second-countable \(\Sigma^0_2\)-compact space.
  Let \(\seq{x_n}_{n=0}^\infty\) be a sequence of elements of \(X.\)
  Consider the \(\Sigma^0_2\) formula \(A(i)\) defined by \[(\exists n_0)(\forall n)(n \geq n_0 \lthen x_n \notin U_i).\]

  We cannot have that \[(\forall x \in X)(\exists i \in I)(A(i) \land x \in U_i).\]
  Indeed, by \(\Sigma^0_2\)-compactness there would be a finite set \(\set{i_0,\dots,i_{p-1}} \subseteq I\) such that \(A(i_0),\dots, A(i_{p-1})\) and \(X = U_{i_0} \cup\cdots\cup U_{i_{p-1}}.\)
  By \Bnd{\Sigma^0_2}, and the definition of \(A(i)\) there is some \(n_0\) such that \[(\forall q < p)(\forall n)(n \geq n_0 \lthen x_n \notin U_{i_q}).\]
  But this is absurd since \(x_{n_0} \in X = U_{i_0} \cup\cdots\cup U_{i_{p-1}}.\)

  It follows that there is some \(x \in X\) with the property that \[(\forall i \in I)(x \in U_i \lthen \lnot A(i)).\]
  But this precisely says that \(x\) is an accumulation point of \(\seq{x_n}_{n=0}^\infty.\)
\end{proof}

\noindent
The next example shows that the hypothesis \Bnd{\Sigma^0_2} is necessary for Proposition~\ref{P:2CompactSeqCompact}.

\begin{example}
  \emph{A finite discrete space which is not sequentially compact from a failure of \Bnd{\Sigma^0_2}.}
  
  By Hirst~\cite{Hirst87}, the failure of \Bnd{\Sigma^0_2} is equivalent to the existence of a bounded sequence \(\seq{x_n}_{n=0}^\infty\) such that \(\set{n \in \N : x_n = x}\) is finite for every \(x.\)
  Let \(m\) be such that \(x_n < m\) for every \(n.\)
  If \(X = \set{0,\dots,m-1}\) is endowed with the discrete topology, then \(\seq{x_n}_{n=0}^\infty\) is a sequence of elements of \(X\) which has no accumulation points.
\end{example}

\section{Discrete Spaces}

A countable second-countable space \((X,\mathcal{U},k)\) is discrete if \(\set{x}\) is a (necessarily basic) open set for every \(x \in X.\)
Again, this notion is not practical over \RCA\ since there may be no effective way to obtain an index for \(\set{x}\) from \(x.\)

\begin{definition}[\RCA]\label{D:EffDiscrete}
  A countable second-countable space \((X,\mathcal{U},k)\) with \(\mathcal{U} = \seq{U_i}_{i \in I}\) is \emph{effectively discrete} if there is a function \(d:X \to I\) such that \(U_{d(x)} = \set{x}\) for every \(x \in X.\)
\end{definition}

\noindent
In \ACA, every discrete space is effectively discrete since the requirements on the index \(d(x)\) are arithmetical.

The following familiar fact is immediate.

\begin{proposition}[\RCA]\label{P:EffDiscreteCompact}
  Let \((X,\mathcal{U},k)\) be a countable second-countable space.
  If \((X,\mathcal{U},k)\) is effectively discrete and compact then \(X\) is finite.
\end{proposition}

\noindent
By increasing the compactness requirement, we can drop the effective discreteness to plain discreteness.

\begin{proposition}[\RCA]\label{P:DiscreteSeqCompact}
  Let \((X,\mathcal{U},k)\) be a countable second-countable space.
  If \((X,\mathcal{U},k)\) is discrete and sequentially compact then \(X\) is finite.
\end{proposition}

\noindent
It follows from Proposition~\ref{P:2CompactSeqCompact} that every countable second-countable space which is discrete and \(\Sigma^0_2\)-compact must be finite under \Bnd{\Sigma^0_2}.

Both Proposition~\ref{P:EffDiscreteCompact} and Proposition~\ref{P:DiscreteSeqCompact} show that every countable second-countable which is discrete and compact must be finite in \ACA.
The next example shows that \ACA\ is necessary to show that every countable second-countable space which is discrete and effectively compact must be finite.

\begin{example}\label{E:DiscreteCompact}
  \emph{An infinite countable second-countable discrete space which is effectively compact from a failure of \(\Sigma^0_1\) comprehension.}

  An enumerable set \(A\) is \emph{hypersimple} if \(A\) is not cofinite and for every sequence \(\seq{d_n}_{n=0}^\infty\) of disjoint finite sets there is an \(n\) such that \(d_n \subseteq A;\) this is a straightforward relativization of the classical notion of hypersimple set due to Post~\cite{Post44}.
  The following claim, which is due to Dekker~\cite{Dekker54} in its unrelativized form, shows that the existence of a hypersimple enumerable set is equivalent to the failure of \(\Sigma^0_1\) comprehension.

  \begin{claim}[\RCA]\label{E:Hypersimple:Existence}
    If \(f:\N\to\N\) is a one-to-one function such that \(f[\N]\) is not a set, then the deficiency set \[D_f = \set{x \in \N : (\exists y)(y > x \land f(y) < f(x))}\] is a hypersimple enumerable set.
  \end{claim}

  \begin{proof}
    It is clear that \(D_f\) is an enumerable set which is not cofinite.
    Suppose \(\seq{d_n}_{n=0}^\infty\) is a sequence of disjoint finite sets such that \(d_n \rem D_f \neq \varnothing\) for every \(n.\)
    To effectively determine whether \(y \in f[\N],\) first find \(n\) such that \(\min\set{f(x) : x \in d_n} \geq y\) and then check whether \(y \in \set{f(x) : x \leq \max d_n}.\)
    This is correct because \(d_n \rem D_f \neq \varnothing.\)
  \end{proof}

  \noindent
  The key property of hypersimple enumerable sets that we will need for this construction is the following.
  
  \begin{claim}[\RCA]\label{E:Hypersimple:Key}
    Let \(A\) be a hypersimple enumerable set.
    For every sequence of finite sets \(\seq{s_n}_{n=0}^\infty\) such that \(s_n \rem A \neq \varnothing\) for every \(n,\) there is a finite set \(b\) such that \(b \cap A = \varnothing\) and \(b \cap s_n \neq \varnothing\) for every \(n.\)
  \end{claim}

  \begin{proof}
    Write \(A = \seq{a_m}_{m=0}^\infty.\)
    Define \(r:\N\times\N\to\N\) by \(r(m,n) = \min(s_n \rem a_m).\)
    Note that \(\lim_{m\to\infty} r(m,n) = \min(s_n \rem A)\) for every \(n.\)
    
    If there is a \(r_0\) such that \(r(m,n) \leq r_0\) for all \(m,n \in \N,\) then the finite set \(b = \set{x \leq r_0 : x \notin A}\) is as required since \(\lim_{m\to\infty} r(m,n) \in b \cap s_n\) for every \(n.\)
    
    Otherwise, recursively define \(\seq{m_i,n_i}_{i=0}^\infty\) such that \[\max(s_{n_i}\rem a_{m_i}) < r(m_{i+1},n_{i+1}) = \min(s_{n_{i+1}} \rem a_{m_{i+1}})\] for every \(i.\)
    Then the sequence \(\seq{d_i}_{i=0}^\infty\) defined by \(d_i = s_{n_i} \rem a_{m_i}\) is a sequence of disjoint finite sets such that \(d_i \rem A \neq \varnothing\) for every \(i,\) which contradicts the fact that \(A\) is hypersimple.
  \end{proof}

  Let \(A = \seq{a_n}_{n=0}^\infty\) be a hypersimple enumerable set.
  Let \(I = \set{\seq{x,y} \in \N^2 : x \leq y}\) and let \(\mathcal{B} = \seq{B_{\seq{x,y}}}_{\seq{x,y} \in I}\) be defined by \[B_{\seq{x,y}} = \set{n \in \N : n = x \lor y \in a_n}.\]
  It is clear that every singleton set is open in \((\N,\mathcal{B}^*).\)
  Indeed, since \(A\) is not cofinite, given \(x\) there is a \(y \geq x\) such that \(y \notin A\) and then \(B_{\seq{x,y}} = \set{x}.\)

  \begin{claim}[\RCA]\label{E:Hypersimple:Compact}
    \((\N,\mathcal{B}^*)\) is compact.
  \end{claim}

  \begin{proof}
    Suppose \(\seq{t_n}_{n=0}^\infty\) is a sequence of elements if \(I^{[<\infty]}\) such that \(\bigcup_{n=0}^\infty B^*_{t_n} = \N.\) 
    We show that there must be some \(n\) such that \(B^*_{t_n}\) is cofinite.

    For each \(n,\) let \(s_n = \set{y : (\exists x \leq y)(\seq{x,y} \in t_n)}.\)
    If there is an \(n\) such that \(s_n \subseteq A,\) then \(B^*_{t_n}\) is cofinite.
    Otherwise, by Claim~\ref{E:Hypersimple:Key}, there must be a finite set \(b\) such that \(b \cap A = \varnothing\) and \(s_n \cap b \neq \varnothing\) for every \(n.\)
    This entails that \(B^*_{t_n} \subseteq \set{0,1,\dots,\max b}\) for every \(n,\) which contradicts the fact that \(\bigcup_{n=0}^\infty B^*_{t_n} = \N.\)
  \end{proof}

  \noindent
  Note that the above proof shows that every open set in \((\N,\mathcal{B}^*)\) is either finite or cofinite.
  Thus the open sets in \((\N,\mathcal{B}^*)\) are precisely the finite or cofinite sets.

  Finally, to show that \((\N,\mathcal{B}^*)\) is effectively compact, we show that \(\mathcal{B}\) has a finite cover relation, which is sufficient by Lemma~\ref{L:SubFiniteCover}.

  \begin{claim}[\RCA]\label{E:Hypersimple:FiniteCover}
    \(\mathcal{B}\) has a finite cover relation.
  \end{claim}

  \begin{proof}
    Of course, the finite cover relation for \(\mathcal{B}\) admits a \(\Pi^0_1\) description.
    It also admits a \(\Sigma^0_1\) description, namely \(\bigcup_{\seq{x,y} \in t} B_{\seq{x,y}} = \N\) if and only if \[(\exists n)(\exists \seq{x_0,y_0} \in t)(y_0 \in a_n \land (\forall x < n)(\exists y)(\seq{x,y} \in t)).\]
    To see that this is correct, suppose that \(\bigcup_{\seq{x,y} \in t} B_{\seq{x,y}} = \N.\)
    There must be a \(\seq{x,y} \in t\) such that \(y \in A.\)
    Let \(n\) be minimal such that there is a \(\seq{x_0,y_0} \in t\) with \(y_0 \in a_n.\)
    Note that \(B_{\seq{x,y}} \subseteq \set{x} \cup [n,\infty)\) for every \(\seq{x,y} \in t.\)
    So for each \(x < n\) there must be a \(y\) such that \(\seq{x,y} \in t.\)
    This shows that \(\bigcup_{\seq{x,y} \in t} B_{\seq{x,y}} = \N\) implies the above \(\Sigma^0_1\) formula; the converse is clear.
    Therefore, \(\Delta^0_1\) comprehension suffices to comprehend the finite cover relation for \(\mathcal{B}.\)
  \end{proof}
\end{example}

\section{Hausdorff Spaces}

A countable second-countable space \((X,\mathcal{U},k)\) is \emph{Hausdorff} if any two points of \(X\) have disjoint (basic) open neighborhoods.
This definition is impractical over \RCA\ since there may be no effective way to get two such neighborhoods.

\begin{definition}[\RCA]
  A countable second-countable space \((X,\mathcal{U},k)\) is \emph{effectively Hausdorff} if there are functions \(h_0,h_1:X \times X \to I\) such that if \(x_0 \neq x_1\) then \(x_0 \in U_{h_0(x_0,x_1)},\) \(x_1 \in U_{h_1(x_0,x_1)}\) and \(U_{h_0(x_0,x_1)} \cap U_{h_1(x_0,x_1)} = \varnothing.\) 
\end{definition}

\noindent
Over \ACA, every countable second-countable Hausdorff space is effectively Hausdorff since the requirements on \(h_0\) and \(h_1\) are arithmetical.
Example~\ref{E:DiscreteOrdered} shows that \ACA\ is necessary to show that every countable second-countable Hausdorff space is effectively Hausdorff.

The following definition gives an elegant characterization of effectively Hausdorff countable second-countable spaces.

\begin{proposition}[\RCA]\label{P:HausdorffDiagonal}
  A countable second-countable space \((X,\mathcal{U},k)\) is effectively Hausdorff if and only if the diagonal is effectively closed in the product space \((X,\mathcal{U},k)\times(X,\mathcal{U},k).\)
\end{proposition}

\begin{proof}
  Suppose that \(A \subseteq I \times I\) is an enumerable set such that \[\set{\seq{x_0,x_1} \in X\times X : x_0 \neq x_1} = {\textstyle\bigcup_{\seq{i_0,i_1} \in A} U_{i_0} \times U_{i_1}}.\]
  Note that we necessarily have \(U_{i_0} \cap U_{i_1} = \varnothing\) for all \(\seq{i_0,i_1} \in A.\)
  Given distinct \(x_0,x_1 \in X,\) we can search through \(A\) for a pair \(\seq{i_0,i_1}\) such that \(\seq{x_0,x_1} \in U_{i_0}\times U_{i_1}.\) 
  Define \(\seq{h_0(x_0,x_1),h_1(x_0,x_1)}\) to be the outcome of this search.

  Conversely, suppose that \(h_0,h_1: X\times X \to I\) witness that \((X,\mathcal{U},k)\) is effectively Hausdorff.
  Then \[\set{\seq{x_0,x_1} \in X \times X : x_0 \neq x_1} = {\textstyle\bigcup_{\seq{i_0,i_1} \in A} U_{i_0}\times U_{i_1}},\] where \(A\) is the enumerable set \[A = \set{\seq{h_0(x_0,x_1),h_1(x_0,x_1)} : \seq{x_0,x_1} \in X \times X \land x_0 \neq x_1}.\qedhere\]
\end{proof}

The familiar fact that compact subspaces of Hausdorff spaces are closed admits the following effective version.

\begin{proposition}[\RCA]\label{P:EffCompactClosed}
  Let \((X,\mathcal{U},k)\) be a countable second-countable effectively Hausdorff space.
  Every effectively compact subspace of \((X,\mathcal{U},k)\) is effectively closed.
\end{proposition}

\begin{proof}
  Write \(\mathcal{U} = \seq{U_i}_{i \in I}\) and let \(h_0,h_1 : X \times X \to I\) witness that \((X,\mathcal{U},k)\) is effectively Hausdorff.
  Suppose \(X_0 \subseteq X\) is relatively effectively compact.
  Let \(X_1 = X \rem X_0\) and let \(C_0 \subseteq I^{[<\infty]}\) be a finite cover relation for \(X_0.\)
  
  First, define a function \(c:X_1\to X_0^{[<\infty]}\) such that \(\set{h_0(x_0,x_1) : x_0 \in c(x_1)} \in C_0\) for every \(x_1 \in X_1.\)
  Such a finite set can always be found since \(X_0\) is relatively compact and \(X_0 \subseteq \bigcup_{x_1 \in X_1} U_{h_0(x_0,x_1)}.\)
  Then, define a function \(f:X_1 \to I\) such that \(x_1 \in U_{f(x_1)} \subseteq \bigcap_{x_0 \in c(x_1)} U_{h_1(x_0,x_1)}\) for every \(x_1 \in X_1.\)
  Note that \(U_{f(x_1)} \cap X_0 = \varnothing\) for every \(x_1 \in X_1\) and so \(X_1 = \bigcup_{x_1 \in X_1} U_{f(x_1)}.\)
\end{proof}

Suppose \((X,\mathcal{U},k)\) is a countable second-countable space which is effectively compact and Hausdorff, but not effectively Hausdorff.
The diagonal is an effectively compact subspace of the product \((X,\mathcal{U},k)\times(X,\mathcal{U},k)\) which is not effectively closed by Proposition~\ref{P:HausdorffDiagonal}.
Example~\ref{E:DiscreteOrdered} shows that such spaces must exist if arithmetic comprehension fails.

Compact Hausdorff spaces always enjoy stronger separation properties.

\begin{definition}[\RCA]
  A countable second-countable space \((X,\mathcal{U},k)\) is \emph{effectively regular} if for every effectively closed subset \(X_0 \subseteq X\) and every point \(x_1 \in X \rem X_0\) there are disjoint effectively open sets \(V_0\) and \(V_1\) such that \(X_0 \subseteq V_0\) and \(x_1 \in V_1.\)
\end{definition}

\begin{definition}[\RCA]
  A countable second-countable space \((X,\mathcal{U},k)\) is \emph{effectively normal} if for every pair of disjoint effectively closed subsets \(X_0, X_1 \subseteq X\) there are disjoint effectively open sets \(V_0\) and \(V_1\) such that \(X_0 \subseteq V_0\) and \(X_1 \subseteq V_1.\)
\end{definition}

\noindent
In \ACA, these correspond to the usual notions of regularity and normality.

The fact that every compact Hausdorff space is regular admits the following effective variant.

\begin{proposition}[\RCA]\label{P:CompactHausdorffRegular}
  Every countable second-countable compact effectively Hausdorff space is effectively regular.
\end{proposition}

\begin{proof}
  Let \(X_0\) be an effectively closed set and let \(x_1 \notin X_0.\)
  Since \(X_0 \subseteq \bigcup_{x_0 \in X_0} U_{h_0(x_0,x_1)}\) there is a finite set \(s \subseteq X_0\) such that \(X_0 \subseteq \bigcup_{x_0 \in s} U_{h_0(x_0,x_1)}.\)
  The neighborhoods \[V_0 = \bigcup_{x_0 \in s} U_{h_0(x_0,x_1)}, \quad V_1 = \bigcap_{x_0 \in s} U_{h_1(x_0,x_1)}\] are as required.
\end{proof}

\noindent
When the space is effectively compact, then we can search for the finite set \(s \subseteq X_0\) and hence we have functions \(r:X \rem X_0 \to I^{[<\infty]}\) and \(s:X \rem X_0 \to I\) such that \(x_1 \in U_{s(x_1)},\) \(X_0 \subseteq \bigcup_{i \in s(x_1)} U_i,\) and \(U_{s(x_1)} \cap \bigcup_{i \in s(x_1)} U_i = \varnothing.\) 
This stronger conclusion allows us to achieve normality.

\begin{proposition}[\RCA]\label{P:EffCompactHausdorffNormal}
  Every countable second-countable effectively compact effectively Hausdorff space is effectively normal.
\end{proposition}

\begin{proof}
  Let \(X_0\) and \(X_1\) be disjoint effectively closed subsets of \(X\) and let \(C_0\) be the finite cover relation for \(X_0.\)
  Let \(s:X_1\to I\) and \(r:X_1\to C_0\) be such that \(x_1 \in U_{s(x_1)}\) and \(U_{s(x_1)} \cap \bigcup_{i \in r(x_1)} U_i = \varnothing\) for each \(x_1 \in X_1.\)
  Since \(X_1 \subseteq \bigcup_{x_1 \in X_1} U_{s(x_1)}\) there is a finite set \(t \subseteq X_1\) such that \(X_1 \subseteq \bigcup_{x_1 \in t} U_{s(x_1)}.\)
  The neighborhoods \[V_0 = \bigcap_{x_1 \in t} \bigcup_{i \in r(x_1)} U_i, \quad V_1 = \bigcup_{x_1 \in t} U_{s(x_1)}\] are as required.
\end{proof}

\noindent
The effective compactness requirement from Proposition~\ref{P:EffCompactHausdorffNormal} can be replaced by \(\Sigma^0_2\)-compactness in the presence of \Ind{\Sigma^0_2}.

\begin{proposition}[\RCA\ + \Ind{\Sigma^0_2}]\label{P:2CompactHausdorffNormal}
  Every countable second-countable \(\Sigma^0_2\)-compact effectively Hausdorff space is effectively normal.
\end{proposition}

\begin{proof}
  Let \(X_0\) and \(X_1\) be disjoint effectively closed subsets of \(X.\)
  By Proposition~\ref{P:CompactHausdorffRegular}, for every \(x_1 \in X_1\) there are \(j \in I\) and a finite set \(c \subseteq I\) such that \(x_1 \in U_j,\) \(X_0 \subseteq \bigcup_{i \in c} U_i,\) and \(U_j \cap \bigcup_{i \in c} U_i = \varnothing.\)
  Let \(A(j)\) be the \(\Sigma^0_2\) formula defined by \[j \in I \land (\exists c \in I^{[<\infty]})(X_0 \subseteq {\textstyle\bigcup_{i \in c} U_i} \land U_j \cap {\textstyle\bigcup_{i \in c} U_i} = \varnothing).\]
  Since \((X,\mathcal{U},k)\) is \(\Sigma^0_2\)-compact, there is a finite set \(\set{j_0,\dots,j_{n-1}} \subseteq I\) such that \(A(j_m)\) holds for every \(m < n\) and \(X_1 \subseteq U_{j_0} \cup\cdots\cup U_{j_{n-1}}.\)
  By \Ind{\Sigma^0_2}, we can find a corresponding sequence \(c_0,\dots,c_{n-1} \in I^{[<\infty]}\) such that \[X_0 \subseteq {\textstyle\bigcup_{i \in c_m} U_i}, \quad U_{j_m} \cap {\textstyle\bigcup_{i \in c_m} U_i} = \varnothing\] for every \(m < n.\)
  Then the neighborhoods \[V_0 = \bigcap_{m=0}^{n-1} \bigcup_{i \in c_m} U_i, \quad V_1 = \bigcup_{m=0}^{n-1} U_{j_m}\] are as required.
\end{proof}

\section{Ordered Spaces}

The order topology on a countable linear order leads to many examples of countable second-countable spaces.
Since all such spaces are Hausdorff, they can be used to answer some of the questions that arose before in our discussion of separation axioms.

\begin{definition}
  Let \((X,{\prec})\) be a countable linear order.
  The associated \emph{ordered space} is the countable second-countable space \((X,\mathcal{U_{\prec}},k_{\prec})\) where \(\mathcal{U}_{\prec} = \seq{(x_0,x_1)}_{\seq{x_0,x_1} \in I},\) \(I = (X\cup\set{\pm\infty})^2,\) and \(k_{\prec}(x,\seq{x_0,x_1},\seq{x'_0,x'_1}) = \seq{\max(x_0,x'_0),\min(x_1,x'_1)}.\)
\end{definition}

A \emph{gap} in \((X,{\prec})\) is a pair \(\seq{x_0,x_1}\) such that \(x_0 \prec x_1\) and \((x_0,x_1) = \varnothing.\)
It turns out that the gap structure of a linear order controls the effective Hausdorffness of the associated ordered space.

\begin{proposition}[\RCA]\label{P:OrderedHausdorff}
  Let \((X,{\prec})\) be a countable linear order.
  The associated ordered space \((X,\mathcal{U}_{\prec},k_{\prec})\) is effectively Hausdorff if and only if the the gaps of \((X,{\prec})\) form a set.
\end{proposition}

\begin{proof}
  Let \(G\) be the set of gaps of \((X,{\prec}).\)
  Define \(h_0,h_1: X \times X \to I\) as follows.
  Assume \(x \prec y.\)
  If \(\seq{x,y} \in G\) then let \(h_0(x,y) = h_1(y,x) = (-\infty,y)\) and \(h_1(x,y) = h_0(y,x) = (x,+\infty).\)
  Otherwise, find the first \(z \in X\) with \(x \prec z \prec y\) and let \(h_0(x,y) = h_1(y,x) = (-\infty,z)\) and \(h_1(x,y) = h_0(y,x) = (z,+\infty).\)
  (For completeness, also let \(h_0(x,x) = h_1(x,x) = (-\infty,+\infty)\) for every \(x \in X.\))
  It is clear that \(h_0,h_1\) witness that \((X,\mathcal{U}_{\prec},k_{\prec})\) is effectively Hausdorff.

  Suppose \(h_0,h_1:X \times X \to I\) witness that \((X,\mathcal{U}_{\prec},k_{\prec})\) is effectively Hausdorff.
  Note that if \(\seq{x_0,x_1}\) is a gap, then \(h_0(x_0,x_1)\) must be an interval with upper endpoint \(x_1\) and \(h_1(x_0,x_1)\) must be an interval with lower endpoint \(x_0.\)
  Conversely, if \(h_0(x_0,x_1)\) is an interval with upper endpoint \(x_1\) and \(h_1(x_0,x_1)\) is an interval with lower endpoint \(x_0,\) then \(\seq{x_0,x_1}\) must be a gap since these intervals do not overlap.
  Thus, \(h_0\) and \(h_1\) provide an effective test to determine whether \(\seq{x_0,x_1}\) is a gap.
\end{proof}

\noindent
The following example shows, in particular, that \ACA\ is necessary to show that every countable second-countable Hausdorff space is effectively Hausdorff.

\begin{example}\label{E:DiscreteOrdered}
  \emph{A discrete ordered space which is not effectively Hausdorff from a failure of \(\Sigma^0_1\) comprehension.}

  Let \(f:\N\to\N\) be a one-to-one function and define the linear order \((\N,{\prec})\) by \(x \prec y\) iff \(f(x) < f(y).\)
  Note that every \(x \in \N\) has an immediate successor \(x^+\) in \((\N,{\prec}),\) and, with the exception of the minimal element, every \(x \in \N\) has an immediate predecessor \(x^-\) in \((\N,{\prec}).\)
  Thus \(\set{x} = (x^-,x^+),\) where \(x^- = -\infty\) in case \(x\) is the minimal element of \((\N,{\prec}).\) 

  Suppose that the gaps of \((\N,{\prec})\) form a set \(G.\)
  Note that the function \(x \mapsto x^+\) is computable from \(G\) since \(x^+\) is the unique element of \((\N,{\prec_f})\) such that \(\seq{x,x^+} \in G.\)
  Let \(x_0\) be the minimal element of \((\N,{\prec})\) and recursively define \(x_{n+1} = x_n^+.\) 
  Then the sequence \(\seq{f(x_n)}_{n=0}^\infty\) is the increasing enumeration of \(f[\N],\) which must therefore be a set.
  
  It follows from Proposition~\ref{P:OrderedHausdorff} that if \(f[\N]\) is not a set, then the ordered space associated to \((\N,{\prec})\) is not effectively Hausdorff.
\end{example}

\noindent
The following simple variation of the preceding example shows that the effective Hausdorffness hypothesis in Proposition~\ref{P:EffCompactClosed} cannot be weakened to plain Hausdorffness.

\begin{example}\label{E:CompactOrdered}
  \emph{An effectively compact ordered space which is not effectively Hausdorff from a failure of \(\Sigma^0_1\) comprehension.}

  Let \(f:\N\to\N\) be a one-to-one function and define the linear order \((\N,{\prec})\) by \(x \prec y\) iff \[(y = 0 \land x \neq 0) \lor (x \neq 0 \land y \neq 0 \land f(x-1) < f(y-1)).\]
  Thus \((\N,{\prec})\) is obtained from Example~\ref{E:DiscreteOrdered} by adding a new maximal element, namely \(0.\)
  Since all basic open neighborhoods of \(0\) are cofinite, the fact that \((\N,\mathcal{U}_{\prec},k_{\prec})\) is \(\Sigma^0_n\)-compact follows immediately from Proposition~\ref{P:CofiniteNCompact} under \Ind{\Sigma^0_n}.
  The fact that \((\N,\mathcal{U}_{\prec},k_{\prec})\) is not effectively Hausdorff when \(f[\N]\) is not a set can be seen in the same way as Example~\ref{E:DiscreteOrdered}.
  Finally, the fact that \((\N,\mathcal{U}_{\prec},k_{\prec})\) has a finite cover relation follows from the following proposition.
\end{example}
 
\begin{proposition}[\RCA]\label{P:OrderedFiniteCover}
  Every countable ordered space has a finite cover relation.
\end{proposition}

\begin{proof}
  Let \((X,{\prec})\) be a countable linear order.
  It is enough to show that the subbase \(\mathcal{B} =
  \seq{B_{\seq{x,i}}}_{\seq{x,i} \in X\times\set{0,1}}\) defined by \[B_{\seq{x,0}} = (-\infty,x), \quad B_{\seq{x,1}} = (x,+\infty)\] has a finite cover relation.
  This is clear since \(t \in (X\times\set{0,1})^{[<\infty]}\) satisfies \(X = \bigcup_{\seq{x,i} \in t} B_{\seq{x,i}}\) if and only if there are \(\seq{x_0,0}, \seq{x_1,1} \in t\) such that \(x_1 \prec x_0.\)
\end{proof}

Recall that a linear order \((X,{\prec})\) is \emph{complete} if for every partition \(X = A^- \cup A^+\) such that \[(\forall x^- \in A^-)(\forall x^+ \in A^+)(x^- \prec x^+),\] either \(A^-\) has a maximal element or \(A^+\) has a minimal element.
Note that if the partition \(A^-,A^+\) witnesses that \((X,{\prec})\) is not complete, then \[\set{(-\infty,x^-) : x_i \in A^-} \cup \set{(x^+,+\infty) : x^+ \in A^+}\] is a basic open cover of \((X,\mathcal{U}_{\prec},k_{\prec})\) which has no finite subcover.
Therefore, every compact ordered space must come from a complete linear order.

Completeness of a countable linear order can be strengthened as follows.

\begin{definition}
  A linear order \((X,{\prec})\) is \emph{\(\Sigma^0_n\)-complete} if for all \(\Sigma^0_n\) formulas \(A^-(x)\) and \(A^+(x)\) such that \[(\forall x^-,x^+ \in X)(A^-(x^-) \land A^+(x^+) \lthen x^- \prec x^+)\] there is an \(x \in X\) such that \[(\forall x^- \in X)(A^-(x^-) \lthen x^- \preceq x) \land (\forall x^+ \in X)(A^+(x^+) \lthen x \preceq x^+).\]
\end{definition}

\noindent
Note that \(\Sigma^0_1\)-completeness is equivalent to plain completeness.
A variation of the argument given above shows that any \(\Sigma^0_n\)-compact ordered space must come from a \(\Sigma^0_n\)-complete linear order.
The following is a partial converse of this fact.

\begin{proposition}[\RCA\ + \Ind{\Sigma^0_{n+1}}]\label{P:CompleteOrdered}
  If \((X,{\prec})\) is a \(\Sigma^0_{n+1}\)-complete linear order, then the associated ordered space \((X,\mathcal{U}_{\prec},k_{\prec})\) is \(\Sigma^0_n\)-compact.
\end{proposition}

\begin{proof}
  We prove the contrapositive.
  Suppose \((X,{\prec})\) is not \(\Sigma^0_n\)-compact.
  Say, \(A(x^-,x^+)\) is a \(\Sigma^0_n\) formula such that \[(\forall x \in X)(\exists x^-,x^+ \in X\cup\set{\pm\infty})(A(x^-,x^+) \land x \in (x^-,x^+))\] but there is no finite sequence \(x^-_0,x^+_0,\ldots,x^-_{\ell-1},x^+_{\ell-1}\) such that \(A(x^-_m,x^+_m)\) for all \(m < \ell\) and \[X = (x^-_0,x^+_0)\cup\cdots\cup(x^-_{\ell-1},x^+_{\ell-1}).\]
  
  Define the \(\Sigma^0_n\) equivalence relation \(\equiv\) on \(X\) by \(x \equiv y\) if and only if there is a finite sequence \(x^-_0,x^+_0,\ldots,x^-_{\ell-1},x^+_{\ell-1}\) such that \(A(x^-_m,x^+_m)\) for all \(m < \ell\) and \[X = (-\infty,\min(x,y))\cup(x^-_0,x^+_0)\cup\cdots\cup(x^-_{\ell-1},x^+_{\ell-1})\cup(\max(x,y),+\infty).\]
  Note that \({\equiv}\)-equivalence classes are convex in the \({\prec}\)-ordering and hence \({\prec}\) naturally induces a linear ordering on the quotient \(X/{\equiv}.\) 
  This quotient is hardly formalizable in second-order arithmetic, but we will still write \([x_0] \prec [x_1]\) as an abbreviation for \(x_0 \prec x_1 \land x_0 \not\equiv x_1.\)

  \begin{claim}
    If \((X,{\prec})\) is \(\Sigma^0_n\)-complete, then \(X/{\equiv}\) is dense-in-itself.
  \end{claim}

  \begin{proof}
    Suppose \(z^-\) and \(z^+\) are such that \([z^-] \prec [z^+].\)
    Then \[(\forall x^-,x^+)(x^- \equiv z^- \land x^+ \equiv z^+ \lthen x^- \prec x^+).\]
    So, by \(\Sigma^0_n\)-completeness, there is a \(z \in X\) such that \[(\forall x^-)(x^- \equiv z^- \lthen x^- \preceq z) \land (\forall x^+)(x^+ \equiv z^+ \lthen z \preceq x^+).\]
    To check that \([z^-] \prec [z] \prec [z^+]\) it suffices to verify that \(z^- \not\equiv z\) and \(z \not\equiv z^+.\)

    Suppose instead that \(z^- \equiv z\) and let \(x^-_0,x^+_0,\ldots,x^-_{\ell-1},x^+_{\ell-1}\) such that \(A(x^-_m,x^+_m)\) for all \(m < \ell\) and \[X = (-\infty,z^-)\cup(x^-_0,x^+_0)\cup\cdots\cup(x^-_{\ell-1},x^+_{\ell-1})\cup(z,+\infty).\]
    Note that \(z'= \max(x^+_0,\ldots,x^+_{\ell-1})\) must be such that \(\seq{z,z'}\) is a gap of \((X,{\prec}).\)
    Otherwise, any point of \((z,z')\) would contradict the defining property of \(z.\)
    We can then find \(x^-,x^+\) such that \(A(x^-,x^+)\) and \(z' \in (x^-,x^+).\)
    But then, \(x^-_0,x^+_0,\ldots,x^-_{\ell-1},x^+_{\ell-1},x^-,x^+\) is a finite sequence which witnesses that \(z^- \equiv z'.\)
    Therefore, \(z'\) contradicts the defining property of \(z.\)
    
    The proof that \(z \not\equiv z^+\) is symmetric.
  \end{proof}

  Now suppose that \((X,{\prec})\) is \(\Sigma^0_n\)-complete.
  Let \(F(x,y^-,y^+,z^-,z^+)\) be the \(\Sigma^0_{n+1}\) formula which says that \(\seq{z^-,z^+}\) is the first pair such that \([y^-]  \prec [z^-] \prec [z^+]  \prec [y^+]\) and either \(x \prec z^-\) or \(z^+ \prec x.\)
  Note that such a pair \(\seq{z^-,z^+}\) is guaranteed to exist when \([y^-] \prec [y^+]\) by density.
 
  Fix an enumeration \(\seq{x_i}_{i=0}^\infty\) of \(X\) such that \(x_0,x_1\) are respectively the minimal and maximal elements of \((X,{\prec})\) (which must exist by \(\Sigma^0_n\)-completeness).
  Let \(Z(i,z^-,z^+)\) be the \(\Sigma^0_{n+1}\) formula which says that there is a sequence \(z_0^-,z_0^+,\dots,z_i^-,z_i^+\) such that \(x_0 = z_0^-, x_1 = z_0^+,\) \(z^- = z_i^-, z^+ = z_i^+,\) and \(F(x_j,z_j^-,z_j^+,z_{j+1}^-,z_{j+1}^+)\) holds for each \(j < i.\)

  By \Ind{\Sigma^0_{n+1}}, for each \(i\) there is a unique pair \(\seq{z^-,z^+}\) such that \(Z(i,z^-,z^+).\)
  Let \(A^-(x^-)\) and \(A^+(x^+)\) be defined by \[A^-(x^-) \liff (\exists i, z^-,z^+)(Z(i,z^-,z^+) \land x^- \preceq z^-),\] and  \[A^+(z^+) \liff (\exists i, z^-,z^+)(Z(i,z^-,z^+) \land z^+ \preceq x^+).\]
  These are \(\Sigma^0_{n+1}\) formulas such that \[(\forall x^-, x^+ \in X)(A^-(x^-) \land A^+(x^+) \lthen x^- \prec x^+)\] and there is no \(x \in X\) such that \[(\forall x^- \in X)(A^-(x^-) \lthen x^- \preceq x) \land (\forall x^+ \in X)(A^+(x^+) \lthen x \preceq x^+). \qedhere\]
\end{proof}

The next example shows that in the case \(n = 1,\) the hypothesis \(\Sigma^0_2\)-completeness cannot be weakened to plain completeness in Proposition~\ref{P:CompleteOrdered}.

\begin{example}\label{E:DiscreteComplete}
  \textit{An infinite complete countable linear order whose associated ordered space is effectively discrete from the existence of a maximal Turing degree \textup(see below\textup).}

  Suppose \(\seq{W_e}_{e=0}^\infty\) is a sequence of enumerable sets such that for every enumerable set \(A\) there is an index \(e\) such that \(A = W_e.\)
  Let us write \(W_e = \seq{w_{e,s}}_{s=0}^\infty\) and assume that \(w_{e,s} \subseteq \set{0,\dots,s-1}\) for every \(s.\)
  Such a two-dimensional array \(\seq{w_{e,s}}_{e,s=0}^\infty\) has maximal Turing degree.
  Indeed, for every set \(X\) there must be indices \(e,f\) such that \(X = W_e\) and \(\N \setminus X = W_f.\)
  Then \(X\) is effectively computable from the two-dimensional array and the indices \(e,f.\)
  Conversely, if there is a set of maximal Turing degree, then one can effectively construct a sequence \(\seq{W_e}_{e=0}^\infty\) of enumerable sets as above.

  This example was greatly inspired by a construction of Watnick~\cite{Watnick84}.
  We construct the linear order \((\N,{\prec})\) by stages, deciding the ordering of finitely many points at a time.
  At each stage, the points of the linear order will be organized into blocks.
  There will be one block \(B_d\) for each dyadic rational \(d\) between \(0\) and \(1,\) inclusive.
  We will always arrange that if \(d < d'\) then the elements of \(B_d\) will precede the elements of \(B_{d'}.\)
  As the construction goes on, the points of the linear order will occasionally shift from one block to another.
  However, each point will eventually stabilize into a fixed block.

  With the exception of blocks \(B_0\) and \(B_1,\) the elements of block \(B_d\) will be indexed consecutively along \(\Z,\) the elements of block \(B_0\) will be indexed consecutively along \(\N\) and the elements of block \(B_1\) will be indexed consecutively along \(-\N = \set{0,-1,-2,\dots}.\)
  As the construction goes on, the blocks will only grow at the ends starting from index \(0.\)
  This way, even as they shift from one block to another, two elements of the same block will always remain at a fixed distance from each other.
  Moreover, as the blocks gradually stabilize, each block \(B_d\) will end up having the same order type as \(\Z,\) except for \(B_0\) and \(B_1\) which will have the same order type as \(\N\) and \(-\N,\) respectively.
  
  For bookkeeping purposes, we prioritize the blocks and the points within each block.
  The \emph{birthday} of the block \(B_d\) is the smallest integer \(n \geq 0\) such that \(d2^n\) is an integer.
  We will write \(D_n\) for the finite set of dyadic rationals in \([0,1]\) with birthday at most \(n,\) i.e., \[D_n = \left\{\frac{0}{2^n},\frac{1}{2^n},\frac{2}{2^n},\dots,\frac{2^n-1}{2^n},\frac{2^n}{2^n}\right\}.\]
  The \emph{level} of the point \(x \in \N\) is the sum of the birthday of the block \(x\) belongs to and the absolute value of the index of \(x\) within that block --- so the element with index \(-3\) in \(B_{3/4}\) has level \(5.\)

  To ensure that the ordered space associated to \((\N,{\prec})\) is effectively discrete, at each stage \(s\) of the construction we will add sufficiently many points so that all points with level at most \(s\) do exist.
  For example, at stage \(5\) we will add new points as necessary to ensure that the block \(B_{3/4}\) does have points with index \(-3,-2,-1,0,1,2,3.\)
  These points may eventually shift to another block, but when that happens new points will be added to fill the vacated spots.
  Since consecutive points in the same block will remain consecutive even when they shift blocks, we can use this property of the construction to effectively compute the two immediate neighbors of any point \(x\) of \((\N,{\prec}).\)
  In turn, these two neighbors give an index for the basic open set \(\set{x},\) as required for effective discreteness.

  For each pair \(W_e,W_f\) of enumerable sets from our universal sequence, we will meet the following requirement.
  \begin{description}
  \item[\((R_{\seq{e,f}})\)]
    If \(W_e \prec W_f\) then there is an \(x\) such that \(W_e \preceq x \preceq W_f.\)
  \end{description}
  The satisfaction of these requirements will guarantee that \((\N,{\prec})\) is a complete linear order. 
  The basic strategy to satisfy requirement \(R_{\seq{e,f}}\) is as follows.
  We will wait for a stage \(s\) such that \(\max w_{e,s} \prec \min w_{f,s}\) and \(\max w_{e,s}\) and \(\min w_{f,s}\) are sufficiently close to each other (defined below).
  At that time, we will shift some points from block to block in such a way that \(\max w_{e,s}\) and \(\min w_{f,s}\) belong to the same block.
  Since \(\max w_{e,s}\) and \(\min w_{f,s}\) will henceforth remain in the same block and at the same finite distance from each other, this will force \(R_{\seq{e,f}}\) to be satisfied.
  
  To determine whether \(\max w_{e,s}\) and \(\min w_{f,s}\) are sufficiently close to each other, we perform the following check.
  If there is at most one block \(B_d\) with birthday at most \(\seq{e,f}\) which intersects the interval \([\max w_{e,s},\min w_{f,s}]_{\prec},\) then and only then can we perform a shift to satisfy requirement \(R_{\seq{e,f}}.\)
  When there is such a block \(B_d,\) then the shift is performed by moving all points from all blocks which intersect the interval \([\max w_{e,s},\min w_{f,s}]_{\prec}\) into \(B_d.\) 
  Note that the points of \(B_d\) are not moved in this process and so this shift does not move any points of level at most \(\seq{e,f}.\)
  When there is no such block, then pick an arbitrary block intersecting \([\max w_{e,s}, \min w_{f,s}]_{\prec}\) to play the role of \(B_d\) as above.
  Again, no points of level at most \(\seq{e,f}\) are moved in this process.

  Of course, when \(\max w_{e,s}\) and \(\min w_{f,s}\) belong to the same block then there is nothing to do for requirement \(R_{\seq{e,f}}.\)
  It follows that there will be at most one shift ever performed to satisfy requirement \(R_{\seq{e,f}}.\)
  By \Ind{\Sigma^0_1}, for every \(n\) there is a stage \(t_n\) after which all shifts for requirements \(R_{\seq{e,f}}\) with \(\seq{e,f} \leq n\) to be performed have been performed.
  Since shifting for requirement \(R_{\seq{e,f}}\) never moves points of level at most \(\seq{e,f},\) the contents of a block \(B_d\) with birthday at most \(n\) are never shifted after stage \(t_n.\)
  Therefore, any points that belong to \(B_d\) at stage \(t_n\) and any points that get into block \(B_d\) after then will forever remain in \(B_d\) and will always keep the same index in \(B_d.\)
  Thus, our promise to have the contents of the block \(B_d\) eventually stabilize will be fulfilled after stage \(t_n.\)

  As observed above, once a shift for requirement \(R_{\seq{e,f}}\) has been performed then this requirement \(R_{\seq{e,f}}\) will automatically be satisfied.
  However, it may very well be that \(\max w_{e,s}\) and \(\min w_{f,s}\) never come sufficiently close to each other at any stage \(s\) and no shifting for requirement \(R_{\seq{e,f}}\) is ever performed.
  It turns out that requirement \(R_{\seq{e,f}}\) will still be satisfied in this scenario.
  Indeed, suppose that at any stage \(s\) there are at least two blocks with birthday at most \(\seq{e,f}\) which intersect \([\max w_{e,s}, \min w_{f,s}]_{\prec}.\)
  There are two ways in which this can occur:
  \begin{enumerate}[(i)]
  \item\label{E:CompleteDiscrete:CaseMid}
    At every stage \(s,\) there is at least one block with birthday at most \(\seq{e,f}\) which is entirely contained in the open interval \((\max w_{e,s}, \min w_{f,s})_{\prec}.\)
  \item\label{E:CompleteDiscrete:CaseEnd}
    From some stage on, there are exactly two blocks with birthday at most \(\seq{e,f}\) which intersect \([\max w_{e,s}, \min w_{f,s}]_{\prec},\) namely the two blocks that contain the endpoints \(\max w_{e,s}\) and \(\min w_{f,s}.\)
 \end{enumerate}
 
 \emph{Ad~\eqref{E:CompleteDiscrete:CaseMid}.}
 Consider the sequence of finite sets \[E_s = \set{ d \in D_{\seq{e,f}} : B_d \subseteq (\max w_{e,s}, \min w_{f,s})_{\prec}}.\]
 Let \(s_0\) be a stage after which blocks with birthday at most \(\seq{e,f}\) are no longer shifted.
 Note that \(\seq{E_s}_{s=s_0}^\infty\) is a decreasing sequence of finite nonempty sets.
 By \(\Ind{\Sigma^0_1},\) the intersection \(E = \bigcap_{s=s_0}^\infty E_s\) must be nonempty.
 If \(d \in E,\) then the point \(x\) of index \(0\) in \(B_d\) satisfies \(\max w_{e,s} \prec x \prec \min w_{f,s}\) for every \(s.\)
 Therefore, this point \(x\) witnesses the conclusion \(W_e \preceq x \preceq W_f\) of requirement \(R_{\seq{e,f}}.\)

 \emph{Ad~\eqref{E:CompleteDiscrete:CaseEnd}.}
 Let \(s_0\) be a stage after which blocks with birthday at most \(\seq{e,f}\) are no longer shifted.
 Further assume that after stage \(s_0,\) there are exactly two blocks with birthday at most \(\seq{e,f}\) which intersect \([\max w_{e,s}, \min w_{f,s}]_{\prec},\) namely the two blocks that contain the endpoints \(\max w_{e,s}\) and \(\min w_{f,s}.\)
 Note that although the points \(\max w_{e,s}\) and \(\min w_{f,s}\) may still change, the two blocks which contain these points must be the same for every \(s \geq s_0\) --- let these blocks be \(B_{d_0}\) and \(B_{d_1},\) respectively.
 Pick a dyadic rational \(d \in (d_0,d_1)\) and let \(x\) be the point with index \(0\) in \(B_d.\)
Since \(\max w_{e,s} \prec x \prec \min w_{f,s}\) for every \(s,\) this point \(x\) witnesses the conclusion \(W_e \preceq x \preceq W_f\) of requirement \(R_{\seq{e,f}}.\)
\end{example}

The existence of a maximal Turing degree is a rather strong requirement, but the above construction works over the minimal \(\omega\)-model \(\mathrm{REC}\) of \RCA.
Hirschfeldt, Shore, and Slaman~\cite{HirschfeldtShoreSlaman09} have shown that the non-existence of a maximal Turing degree is equivalent to the Atomic Model Theorem with Subenumerable Types (\(\mathsf{AST}\)).
It would be interesting to know whether \(\mathsf{AST}\) is strong enough to show that the ordered space associated to a complete linear ordering is compact.

\section{The Tychonoff Theorem}

It is straightfoward to see that the product of sequentially compact spaces is sequentially compact.

\begin{proposition}[\RCA]
  If \((X,\mathcal{U},k)\) and \((Y,\mathcal{V},\ell)\) are countable second-countable sequentially compact spaces, then the product \((X,\mathcal{U},k)\times(Y,\mathcal{V},\ell)\) is sequentially compact too.
\end{proposition}

\begin{proof}
  Let \(\seq{\seq{x_n,y_n}}_{n=0}^\infty\) be a sequence of elements of \(X \times Y.\)
  By hypothesis, the sequence \(\seq{x_n}_{n=0}^\infty\) has an accumulation point \(x \in X.\)
  Recursively construct a subsequence \(\seq{x_{n(m)}}_{m=0}^\infty\) with the property that \(\set{m: x_{n(m)} \in U_i}\) is cofinite for each basic neighborhood \(U_i\) of \(x.\)
  By hypothesis, the sequence \(\seq{y_{n(m)}}_{m=0}^\infty\) has an accumulation point \(y.\)
  Then, \(\seq{x,y}\) is an accumulation point of \(\seq{\seq{x_n,y_n}}_{n=0}^\infty.\)
\end{proof}

\noindent
In \ACA, this shows that the product of two countable second-countable compact spaces is compact.
However, \ACA\ is necessary for this approach to work as Example~\ref{E:DiscreteCompact} shows.
In fact, all approaches to prove the Tychonoff Theorem based on accumulation points or convergence, such as Tychonoff's original proof~\cite{Tychonoff35}, are bound to fail below \ACA\ for this reason.

A more elementary approach to the Tychonoff Theorem due to Loeb~\cite{Loeb65} turns out to be much more successful in \RCA.

\begin{lemma}[\RCA]\label{L:CompactLoeb}
  Let \((X,\mathcal{U},k)\) and \((Y,\mathcal{V},\ell)\) be countable second-countable spaces such that \(\mathcal{U} = \seq{U_i}_{i \in I}\) is compact and \(\mathcal{V} = \seq{V_j}_{j \in J}\) has a finite cover relation.
  If \(A \subseteq I \times J\) is an enumerable set such that \(X \times Y \neq \bigcup_{\seq{i,j} \in a} U_i \times V_j\) for every finite set \(a \subseteq A,\) then there is an \(x \in X\) such that \(Y \neq \bigcup_{j \in b} V_j\) for every finite subset \(b\) of the enumerable set \[B = \set{j \in J : (\exists i \in I)(x \in U_i \land \seq{i,j} \in A)}.\]
\end{lemma}

\begin{proof}
  Suppose, for the sake of contradiction, that for every \(x \in X\) there is a finite set \(a \subseteq \set{\seq{i,j} \in A : x \in U_i}\) such that \(Y = \bigcup_{\seq{i,j} \in a} V_j.\)
  Since \(\mathcal{V}\) has a finite cover relation, we can effectively search for such \(a\) and hence we have a function \(f:X \to A^{[<\infty]}\) which witnesses this.
  
  Let \(h:X \to I\) be such that \(x \in U_{h(x)} \subseteq \bigcap_{\seq{i,j} \in f(x)} U_i.\)
  Since \(\mathcal{U}\) is compact, there are \(x_0,\dots,x_{n-1} \in X\) such that \(X = U_{h(x_0)} \cup\cdots\cup U_{h(x_{n-1})}.\)
  Let \(a = f(x_0) \cup\cdots\cup f(x_{n-1}).\) 
  Then \(a\) is a finite subset of \(A\) such that \(X \times Y = \bigcup_{\seq{i,j} \in a} U_i \times V_j.\)
\end{proof}

\noindent
It follows immediately that \RCA\ proves the following weak version of the Tychonoff Theorem.

\begin{proposition}[\RCA]\label{P:HalfTychonoff}
  If \((X,\mathcal{U},k)\) and \((Y,\mathcal{V},\ell)\) are countable second-countable spaces such that \((X,\mathcal{U},k)\) is compact and \((Y,\mathcal{V},\ell)\) is effectively compact, then the product space \((X,\mathcal{U},k)\times(Y,\mathcal{V},\ell)\) is compact.
\end{proposition}

\noindent
Observing that the product of two countable second-countable spaces with finite cover relations has a finite cover relation, we obtain the Tychonoff Theorem for countable second-countable effectively compact spaces.

\begin{proposition}[\RCA]\label{P:WeakTychonoff}
  If \((X,\mathcal{U},k)\) and \((Y,\mathcal{V},\ell)\) are countable second-countable effectively compact spaces, then the product space \((X,\mathcal{U},k)\times(Y,\mathcal{V},\ell)\) is effectively compact.
\end{proposition}

For \(\Sigma^0_2\)-compact spaces, we have the following analogue of Lemma~\ref{L:CompactLoeb}.

\begin{lemma}[\RCA\ + \Bnd{\Sigma^0_2}]\label{L:2CompactLoeb}
  Let \((X,\mathcal{U},k)\) and \((Y,\mathcal{V},\ell)\) be countable second-countable spaces such that \(\mathcal{U} = \seq{U_i}_{i \in I}\) is \(\Sigma^0_2\)-compact and \(\mathcal{V} = \seq{V_j}_{j \in J}.\)
  If \(A \subseteq I \times J\) is an enumerable set such that \(X \times Y \neq \bigcup_{\seq{i,j} \in a} U_i \times V_j\) for every finite set \(a \subseteq A,\) then there is an \(x \in X\) such that \(Y \neq \bigcup_{j \in b} V_j\) for every finite subset \(b\) of the enumerable set \[B = \set{j \in J : (\exists i \in I)(x \in U_i \land \seq{i,j} \in A)}.\]
\end{lemma}

\begin{proof}
  Suppose, for the sake of contradiction, that for every \(x \in X\) there is a finite set \(a \subseteq \set{\seq{i,j} \in A : x \in U_i}\) such that \(Y = \bigcup_{\seq{i,j} \in a} V_j.\)
  Consider the \(\Sigma^0_2\) formula \(C(i_0)\) defined by \[i_0 \in I \land (\exists a \in A^{[<\infty]})(U_{i_0} \subseteq {\textstyle\bigcap_{\seq{i,j} \in a} U_i} \land Y = {\textstyle\bigcup_{\seq{i,j} \in a} V_j}).\]
  Since \((X,\mathcal{U},k)\) is \(\Sigma^0_2\)-compact, there is a finite set \(c \subseteq I\) such that \((\forall i \in c)C(i)\) and \(X = \bigcup_{i \in c} U_i.\) 
  By \Bnd{\Sigma^0_2}, there is a finite set \(a \subseteq A\) such that for every \(i \in c\) there is a finite set \(a_i \subseteq a\) such that \[U_i \subseteq {\textstyle\bigcap_{\seq{i,j} \in a_i} U_i} \land Y = {\textstyle\bigcup_{\seq{i,j} \in a_i} V_j}.\]
  But then \(X \times Y = \bigcup_{\seq{i,j} \in a} U_i \times V_j,\) contrary to the hypothesis.
\end{proof}

\noindent
Consequently, we obtain another weak form of the Tychonoff Theorem.

\begin{proposition}[\RCA\ + \Bnd{\Sigma^0_2}]\label{P:2Tychonoff}
  If \((X,\mathcal{U},k)\) and \((Y,\mathcal{V},\ell)\) are countable second-countable spaces such that \((X,\mathcal{U},k)\) is \(\Sigma^0_2\)-compact and \((Y,\mathcal{V},\ell)\) is compact, then the product space \((X,\mathcal{U},k)\times(Y,\mathcal{V},\ell)\) is compact. 
\end{proposition}

\noindent
In fact, straightforward adaptation of Lemma~\ref{L:2CompactLoeb} yields the following.

\begin{proposition}[\RCA\ + \Bnd{\Sigma^0_n}, \(n \geq 2\)]\label{P:NTychonoff}
  If \((X,\mathcal{U},k)\) and \((Y,\mathcal{V},\ell)\) are countable second-countable \(\Sigma^0_n\)-compact spaces, then the product space \((X,\mathcal{U},k)\times(Y,\mathcal{V},\ell)\) is \(\Sigma^0_n\)-compact.
\end{proposition}

Thus, modulo some inductive assumptions, the Tychonoff Theorem is provable in \RCA\ for all of the strong forms of compactness considered in this paper.
However, the Tychonoff Theorem for plain compactness is not provable in \RCA\ (with any amount of induction) as the following example will show.

\begin{example}\label{E:Tychonoff}
  \textit{A failure of the Tychonoff Theorem for countable second-countable compact spaces from the existence of a \(\Delta^0_2\)-definable \(1\)-generic function \textup(defined below\textup).}

  Assume \Bnd{\Sigma^0_2}.
  Suppose \(\seq{f_n}_{n=0}^\infty\) is a sequence of functions \(f_n:\N\to\set{0,1}\) such that \(\lim_{n\to\infty} f_n(x)\) exists for every \(x.\)
  Note that the function defined by this limit need not exist in \RCA, but \Bnd{\Sigma^0_2} guarantees that all of its initial segments \(\lim_{n\to\infty} f_n\res\ell\) do exist.
  Suppose further that the \(\Delta^0_2\)-definable function defined by this limit is \emph{\(1\)-generic}: for every enumerable set \(D \subseteq \set{0,1}^{<\infty}\) there is an \(\ell\) such that either \(\lim_{n\to\infty} f_n\res\ell \in D,\) or else \(\lim_{n\to\infty} f_n\res\ell\) has no extension in \(D\) at all.
  (This is a relativization of the notion of \(1\)-generic set due to Jockusch and Posner~\cite{JockuschPosner78}.)

  \begin{claim}[\RCA\ + \Bnd{\Sigma^0_2}]\label{E:Tychonoff:Key}
    For \(i \in \set{0,1}\) and every sequence \(\seq{s_m}_{m=0}^\infty\) of finite sets, if it is the case that for every \(m\) there is an \(x \in s_m\) with \(\lim_{n\to\infty} f_n(x) = i,\) then there is a finite set \(b\) such that \(\lim_{n\to\infty} f_n(x) = i\) for every \(x \in b,\) and \(s_m \cap b \neq \varnothing\) for every \(m.\)
  \end{claim}

  \begin{proof}
    Let \(D_i\) be the enumerable set of all \(t \in \set{0,1}^{<\infty}\) such that \(s_m \subseteq t^{-1}(1-i)\) for some \(m.\)
    Since \(\lim_{n\to\infty} f_n\) is \(1\)-generic, there is an \(\ell\) such that either \(\lim_{n\to\infty} f_n\res\ell \in D_i\) or else \(\lim_{n\to\infty} f_n\res\ell\) has no extensions in \(D_i\) at all.
    In the first case, when \(\lim_{n\to\infty} f_n\res\ell \in D_i,\) we have that for some \(m,\) \(\lim_{n\to\infty} f(x) = 1-i\) for every \(x \in s_m.\)
    In the second case, when \(\lim_{n\to\infty} f_n\res\ell\) has no extensions in \(D_i,\) it must be that for every \(m\) there is an \(x \in s_m \cap \set{0,\dots,\ell-1}\) such that \(\lim_{n\to\infty} f_n(x) = i.\)
    The finite set \(b = \set{x < \ell : \lim_{n\to\infty} f_n(x) = i}\) (which exists by \Bnd{\Sigma^0_2}) is then as required.
  \end{proof}

  Let \(I = \set{\seq{x,y} \in \N^2 : x \geq y}.\)
  For \(i \in \set{0,1},\) define \(\mathcal{B}_i = \seq{B_{\seq{x,y}}^i}_{\seq{x,y} \in I}\) by \[B^i_{\seq{x,y}} = \set{n \in \N : n = y \lor f_n(x) = i}.\]
  We will show that the spaces \((\N,\mathcal{B}_0^*)\) and \((\N,\mathcal{B}_1^*)\) are both compact, but that the product \((\N,\mathcal{B}_0^*)\times(\N,\mathcal{B}_1^*)\) is not.

  \begin{claim}[\RCA\ + \Bnd{\Sigma^0_2}]\label{E:Tychonoff:Compact}
    The spaces \((\N,\mathcal{B}_0^*)\) and \((\N,\mathcal{B}_1^*)\) are both compact.
  \end{claim}

  \begin{proof}
    Suppose \(\seq{t_m}_{m=0}^\infty\) is a sequence of elements if \(I^{[<\infty]}\) such that \(\bigcup_{m=0}^\infty B^{i*}_{t_m} = \N.\) 
    We show that there must be some \(m\) such that \(B^{i*}_{t_m}\) is cofinite.

    For each \(m,\) let \(s_m = \set{x : (\exists y \leq x)(\seq{x,y} \in t_m)}.\)
    If there is an \(m\) such that \(\lim_{n\to\infty} f_n(x) = i\) for every \(x \in s_n,\) then \(B^{i*}_{t_n}\) is cofinite by \Bnd{\Sigma^0_2}.
    Otherwise, by Claim~\ref{E:Tychonoff:Key}, there must be a finite set \(b \subseteq \set{x : \lim_{n\to\infty} f_n(x) = 1-i}\) such that \(s_m \cap b \neq \varnothing\) for every \(m.\)
    This entails that \[B^{i*}_{t_n} \subseteq \set{0,\dots,\max b} \cup {\textstyle\bigcup_{x \in b} \set{n : f_n(x) = i}}\] for every \(n.\) 
    Since this superset is finite by \Bnd{\Sigma^0_2}, this contradicts the fact that \(\bigcup_{n=0}^\infty B^*_{t_n} = \N.\)
  \end{proof}

  \begin{claim}[\RCA\ + \Bnd{\Sigma^0_2}]
    The product space \((\N,\mathcal{B}_0^*)\times(\N,\mathcal{B}_1^*)\) is not compact.
  \end{claim}

  \begin{proof}
    Consider the sequence of basic open sets \(\seq{B^0_{\seq{x,y_0}}\times B^1_{\seq{x,y_1}}}_{\seq{x,y_0,y_1} \in A}\) where \(A = \set{\seq{x,y_0,y_1} : x \geq \max(y_0,y_1)}.\)
    We clearly have \[\N\times\N = \bigcup_{\seq{x,y_0,y_1} \in A} B^0_{\seq{x,y_0}} \times B^1_{\seq{x,y_1}}.\]
    We claim that \(\N\times\N \neq \bigcup_{\seq{x,y_0,y_1} \in a} B^0_{\seq{x,y_1}} \times B^1_{\seq{x,y_1}}\) for every finite set \(a \subseteq A.\)
    By \Bnd{\Sigma^0_2}, we can write \(a = a_0 \cup a_1\) where \[a_i = \set{\seq{x,y_0,y_1} \in a : \lim_{n\to\infty} f_n(x) = 1 - i}.\]
    Notice that \(B^i_{\seq{x,y_i}}\) is finite for every \(\seq{x,y_0,y_1} \in a_i.\)
    By \Bnd{\Sigma^0_2} again, we can find \(z\) such that \[\bigcup_{\seq{x,y_0,y_1} \in a_0} B^0_{\seq{x,y_0}} \cup \bigcup_{\seq{x,y_0,y_1} \in a_1} B^1_{\seq{x,y_1}} \subseteq \set{0,\dots,z-1}.\]
    Clearly, \(\seq{z,z} \notin \bigcup_{\seq{x,y_0,y_1} \in a} B^0_{\seq{x,y_0}} \times B^1_{\seq{x,y_1}}.\)
  \end{proof}

  The spaces \((\N,\mathcal{B}_0^*)\) and \((\N,\mathcal{B}_1^*)\) are not particularly nice, but a small variation of the above gives two discrete spaces with the same properties.
  For this variant, use the subbases \(\mathcal{C}_i = \seq{C^i_{\seq{w,x,y}}}_{\seq{w,x,y} \in J}\) where \[C^i_{\seq{w,x,y}} = \set{n \in \N : n = y \lor (n \geq w \land f_n(x) = i)}\] and \(J = \set{\seq{w,x,y} : w \geq x \geq y}.\)
\end{example}

The standard construction of a \(\mathbf{0}'\)-computable \(1\)-generic Turing degree (essentially due to Kleene and Post~\cite{KleenePost54}) shows that there is a \(\Delta^0_2\)-definable \(1\)-generic function over \(\mathrm{REC},\) the minimal \(\omega\)-model of \RCA.
In fact, there is an \(\omega\)-model of \WKL\ over which there is a \(\Delta^0_2\)-definable \(1\)-generic function.
Indeed, there is an \(\omega\)-model of \WKL\ all of whose sets are \(\mathbf{d}\)-computable for some low Turing degree \(\mathbf{d}.\) 
The standard construction of a \(1\)-generic Turing degree relativized to \(\mathbf{d}\) gives a \(1\)-generic function which is computable from \(\mathbf{0}' = \mathbf{d}'.\)
Thus, \WKL\ does not suffice to prove the Tychonoff Theorem for countable second-countable compact spaces.

One of the weakest known principles which prohibits the existence of a \(\Delta^0_2\)-definable \(1\)-generic function is the Stable Ramsey Theorem for pairs (\(\mathsf{SRT}^2_2\)) of Cholak, Jockusch, and Slaman~\cite{CholakJockuschSlaman01}.
It would be interesting to know whether \(\mathsf{SRT}^2_2\) implies that the product of two compact countable second-countable spaces is compact.

\section{The Alexander Subbase Theorem}

There is yet another common approach to prove the Tychonoff Theorem, which goes through the Alexander Subbase Theorem~\cite{Alexander39}.
This approach turns out to be less efficient than Loeb's approach, but the analyis of the method is interesting.

\begin{definition}[\RCA]
  A sequence \(\seq{B_i}_{i \in I}\) of subsets of \(X\) is \emph{compact} if for every enumerable set \(A \subseteq I\) such that \(X = \bigcup_{i \in A} B_i\) there is a finite set \(a \subseteq A\) such that \(X = \bigcup_{i \in a} B_i.\)
  We say that \(\seq{B_i}_{i \in I}\) is \emph{effectively compact} if moreover \(\seq{B_i}_{i \in I}\) has a finite cover relation.
\end{definition}

\noindent
Thus a countable second-countable space \((X,\mathcal{U},k)\) is (effectively) compact if and only if \(\mathcal{U}\) is (effectively) compact as a sequence of subsets of \(X.\)

The Alexander Subbase Theorem says that every space with a compact subbase is compact; this fact is provable in \ACA.

\begin{proposition}[\ACA]\label{P:Alexander}
  If \(\mathcal{B}\) is a compact sequence of subsets of \(X,\) then the countable second-countable space \((X,\mathcal{B}^*)\) is compact.
\end{proposition}

\begin{proof}
  Let \((X,\mathcal{B}^*)\) be a countable second-countable space where \(\mathcal{B} = \seq{B_i}_{i \in I}\) is a compact subbase.

  Let \(A \subseteq I^{[<\infty]}\) be such that \(X \neq \bigcup_{s \in a} B^*_s\) for every finite set \(a \subseteq A.\)
  We will show that \(X \neq \bigcup_{s \in A} B^*_s.\)

  Note that not having a finite subset which is a cover is an arithmetical property of finite character.
  By Tukey's Lemma, which was analyzed by Dzhafarov and Mummert~\cite{DzhafarovMummertXX}, we can assume that for every \(t \in I^{[<\infty]} \rem A\) there is a finite set \(a \subseteq A\) such that \(X = B^*_t \cup \bigcup_{s \in a} B^*_s.\) 
  
  Let \(A_1 = \set{s \in A : |s| = 1}\) and note that \(X \neq \bigcup_{s \in A_1} B^*_s\) since \(B^*_s \in \mathcal{B}\) for each \(s \in A_1.\)
  Pick \(x \in X \rem \bigcup_{s \in A_1} B^*_s;\) we claim that \(x \notin \bigcup_{s \in A} B^*_s.\)
  
  Suppose instead that there is a \(t \in A\) with \(x \in B^*_t.\)
  Write \(t = \set{t_0,\dots,t_{\ell-1}}\) so that \(B^*_t = B_{t_0} \cap\cdots\cap B_{t_{\ell-1}}.\)
  Observe that \(\set{t_k} \notin A\) for each \(k < \ell.\)
  By maximality of \(A,\) for each \(k < \ell\) there is a finite set \(a_k \subseteq A\) such that \(X = B_{t_k} \cup \bigcup_{s \in a_k} B^*_s.\)
  Then \(a = a_0 \cup\cdots\cup a_{\ell-1}\) is a finite subset of \(A\) such that \(X = B^*_t \cup \bigcup_{s \in a} B^*_s,\) which is impossible.
\end{proof}

\noindent
The next example reverses the previous proposition.

\begin{example}\label{E:DiscreteSubcompact}
  \emph{An infinite effectively discrete space with a compact subbase from a failure of \(\Delta^0_2\) comprehension.}

  Let \(\seq{f_n}_{n=0}^\infty\) be a sequence of functions \(f_n:\N\to\set{0,1}\) such that \(\lim_{x\to\infty} f_n(x)\) exists for every \(n.\)
  Define \(\mathcal{B} = \seq{B_s}_{s \in \set{0,1}^{<\infty}}\) by \[B_s = \set{x \in \N : x = |s| \lor (\exists n < |s|)(f_n(x) \neq s_n)}.\]
  Note that \((\N,\mathcal{B}^*)\) is effectively discrete since \(B^*_{\set{0,1}^n} = \set{n}\) for each \(n.\)

  Suppose that \(\mathcal{B}\) is not compact; say \(A \subseteq \set{0,1}^{<\infty}\) is an enumerable set such that \(\N = \bigcup_{s \in A} B_s\) but \(\N \neq \bigcup_{s \in a} B_s\) for every finite set \(a \subseteq A.\)
  Note that if \(s \in A\) then we must have that \(\lim_{x\to\infty} f_n(x) = s_n\) for every \(n < |s|,\) otherwise there would be a \(s \in A\) such that \(B_s\) is cofinite.
  Since \(A\) is infinite, it must contain sequences of arbitrarily long length, thus \[\lim_{x\to\infty} f_n(x) = d \liff (\exists s \in A)(n < |s| \land s_n = d),\] which allows us to comprehend the function \(n \mapsto \lim_{x\to\infty} f_n(x).\)

  Therefore, if \(\seq{f_n}_{n=0}^\infty\) witnesses the failure of \(\Delta^0_2\) comprehension, then \((\N,\mathcal{B}^*)\) is an infinite effectively discrete space with a compact subbase.
\end{example}

For effectively compact subbases, the Weak K{\"o}nig Lemma suffices to prove the Alexander Subbase Theorem.

\begin{proposition}[\WKL]\label{P:WeakAlexander}
  If \(\mathcal{B}\) is an effectively compact sequence of subsets of \(X,\) then the countable second-countable space \((X,\mathcal{B}^*)\) is effectively compact.
\end{proposition}

\begin{proof}
  Write \(\mathcal{B} = \seq{B_i}_{i \in I}\) and let \(C \subseteq I^{[<\infty]}\) be the finite cover relation for \(\mathcal{B}.\)
  We know from Lemma~\ref{L:SubFiniteCover} that \(\mathcal{B}^*\) has a finite cover relation, so it suffices to show that \(\mathcal{B}^*\) is compact.

  Let \(\seq{s_n}_{n=0}^\infty\) be a sequence of elements of \(I^{[<\infty]}\) such that \(X = \bigcup_{n=0}^\infty B^*_{s_n}.\)
  Let \(T\) be the bounded tree of all sequences \(\seq{i_0,\dots,i_{\ell-1}}\) such that \(i_n \in s_n\) for each \(n < \ell\) and such that \(\set{i_0,\dots,i_{\ell-1}} \notin C.\)

  Note that \(T\) does not have any infinite branches.
  Indeed, if \(\seq{i_n}_{n=0}^\infty\) were to enumerate an infinite branch through \(T,\) then we would have \(\bigcup_{n=0}^\infty B_{i_n} = X\) since \(B^*_{s_n} \subseteq B_{i_n}\) for every \(n,\) but \(\bigcup_{n=0}^{\ell-1} B_{i_n} \neq X\) for every \(\ell\) since \(\set{i_0,\dots,i_{\ell-1}} \notin C.\)
 
  It follows from the Weak K{\"o}nig Lemma that \(T\) must be finite. 
  Let \(\ell\) be the height of \(T,\) then
  \[\bigcup_{n=0}^{\ell-1} B^*_{s_n} = \bigcup_{n=0}^{\ell-1} \bigcap_{i \in s_n} B_i = \bigcap_{\seq{i_0,\dots,i_{\ell-1}} \in s_0\times\cdots\times s_{\ell-1}} \bigcup_{n=0}^{\ell-1} B_{i_n} = X,\] with the usual convention that empty intersections equal \(X.\)
\end{proof}

\noindent
The next example reverses the previous proposition.

\begin{example}\label{E:DiscreteEffSubcompact}
  \emph{An infinite effectively discrete space with an effectively compact subbase from a failure of the Weak K{\"o}nig Lemma.}

  Let \(T\) be a subtree of \(\set{0,1}^{<\infty}\) which has no infinite branches and let \(X\) be the set of all dead-ends of \(T.\)
  Let \(\mathcal{B} = \seq{B_t}_{t \in T}\) be defined by \[B_t = \set{x \in X : t \nsubseteq x}.\]
  This is a compact sequence with a finite cover relation.

  The countable second-countable space \((X,\mathcal{B}^*)\) is effectively discrete since \(\set{x} = \bigcap_{t \in b_x} B_t\) where \(b_x = \set{(x \res i)\cat(1-x(i)) : i < |x|} \cap T.\)
  Thus \((X,\mathcal{B}^*)\) is compact if and only if \(X\) (and hence \(T\)) is finite by Proposition~\ref{P:EffDiscreteCompact}.
  Therefore, if \(T\) witnesses the failure of the Weak K{\"o}nig Lemma, then \((X,\mathcal{B}^*)\) is an infinite effectively discrete space with an effectively compact subbase.
\end{example}

\section{Conclusion}

In this paper, we considered compactness, discreteness, and Hausdorffness for countable second-countable spaces in weak subsystems of second-order arithmetic.
We showed that these properties behave as expected in the subsystem \ACA.
In each case, we also found effective variants of these properties and we showed that these are well behaved even in the base subsystem \RCA.

Examples~\ref{E:Compact}, \ref{E:DiscreteCompact}, \ref{E:DiscreteOrdered}, and~\ref{E:CompactOrdered} exist in any model of \RCA\ where arithmetic comprehension fails.
The following table summarizes their key properties:
\begin{center}
  \begin{tabular}{|l|c|c|c|c|}
    \hline
    & \ref{E:Compact} 
    & \ref{E:DiscreteCompact} 
    & \ref{E:DiscreteOrdered}
    & \ref{E:CompactOrdered} \\\hline
    Compact/Effectively & Yes/No & Yes/Yes & No/No & Yes/Yes \\
    Discrete/Effectively & No/No & Yes/No & Yes/No & No/No \\
    Hausdorff/Effectively & No/No & Yes/No & Yes/No & Yes/No \\\hline
  \end{tabular}
\end{center}
It follows that the following statements are all equivalent to arithmetic comprehension over \RCA.
\begin{itemize}
\item Every compact countable second-countable space is effectively compact. 
\item Every discrete (effectively compact) countable second-countable space is effectively discrete.
\item Every Hausdorff (effectively compact) countable second-countable space is effectively Hausdorff.
\end{itemize}
Thus \ACA\ is both necessary and sufficient for the good behavior of compactness, discreteness, and Hausdorffness for countable second-countable spaces.

Example~\ref{E:DiscreteCompact} also shows that the following statements are equivalent to arithmetic comprehension over \RCA.
\begin{itemize}
\item Every discrete (effectively) compact countable second-countable space is finite.
\item Every (effectively) compact countable second-countable space is sequentially compact.
\end{itemize}

Examples~\ref{E:DiscreteComplete} and~\ref{E:Tychonoff} exist under stronger hypotheses than the mere failure of arithmetic comprehension.
These two examples respectively show that:
\begin{itemize}
\item \RCA\ does not prove that the order topology of a countable complete linear ordering is compact.
\item \WKL\ does not prove that the product of two compact countable second-countable spaces is compact.
\end{itemize}
However, the precise strength of these two statements is unknown.
Considering the strength of the hypotheses needed for these two examples, we are led to the following two questions.

\begin{question}
  Does the principle \(\mathsf{AST}\) of Hirschfeldt, Shore, and Slaman~\cite{HirschfeldtShoreSlaman09} imply that the order topology of a countable complete linear ordering is compact?
\end{question}

\begin{question}
  Does the principle \(\mathsf{SRT}^2_2\) of Cholak, Jockusch, and Slaman~\cite{CholakJockuschSlaman01} imply that the product of two compact countable second-countable spaces is compact?
\end{question}

\bibliographystyle{amsplain}
\bibliography{tychonoff}

\providecommand{\bysame}{\leavevmode\hbox to3em{\hrulefill}\thinspace}
\providecommand{\MR}{\relax\ifhmode\unskip\space\fi MR }
\providecommand{\MRhref}[2]{%
  \href{http://www.ams.org/mathscinet-getitem?mr=#1}{#2}
}
\providecommand{\href}[2]{#2}
\begin{thebibliography}{10}

\bibitem{Alexander39}
J.~W. Alexander, \emph{Ordered sets, complexes, and the problem of
  bicompactification}, Proc. Nat. Acad. Sci. U.S.A. \textbf{22} (1939),
  296--298.

\bibitem{CholakJockuschSlaman01}
P.~A. Cholak, C.~G. Jockusch, Jr., and T.~A. Slaman, \emph{On the strength of
  {R}amsey's theorem for pairs}, J. Symbolic Logic \textbf{66} (2001), no.~1,
  1--55. \MR{1825173 (2002c:03094)}

\bibitem{Dekker54}
J.~C.~E. Dekker, \emph{A theorem on hypersimple sets}, Proc. Amer. Math. Soc.
  \textbf{5} (1954), 791--796. \MR{0063995 (16,209b)}

\bibitem{DzhafarovMummertXX}
D.~D. Dzhafarov and C.~Mummert, \emph{Reverse mathematics and properties of
  finite character}, preprint (2011), 15, arXiv:math-lo/1109.3378.

\bibitem{HirschfeldtShoreSlaman09}
D.~R. Hirschfeldt, R.~A. Shore, and T.~A. Slaman, \emph{The atomic model
  theorem and type omitting}, Trans. Amer. Math. Soc. \textbf{361} (2009),
  no.~11, 5805--5837. \MR{2529915 (2010m:03019)}

\bibitem{Hirst87}
J.~L. Hirst, \emph{Combinatorics in subsystems of second order arithmetic},
  Ph.D. thesis, The Pennsylvania State University, 1987, p.~153.

\bibitem{Hunter08}
J.~Hunter, \emph{Higher-order reverse topology}, Ph.D. thesis, The University
  of Wisconsin - Madison, 2008, p.~97.

\bibitem{JockuschPosner78}
C.~G. Jockusch, Jr. and D.~B. Posner, \emph{Double jumps of minimal degrees},
  J. Symbolic Logic \textbf{43} (1978), no.~4, 715--724. \MR{518677
  (80d:03042)}

\bibitem{KleenePost54}
S.~C. Kleene and E.~L. Post, \emph{The upper semi-lattice of degrees of
  recursive unsolvability}, Ann. of Math. (2) \textbf{59} (1954), 379--407.
  \MR{0061078 (15,772a)}

\bibitem{Loeb65}
P.~A. Loeb, \emph{A new proof of the {T}ychonoff theorem}, Amer. Math. Monthly
  \textbf{72} (1965), 711--717. \MR{0190896 (32 \#8306)}

\bibitem{Mummert06}
C.~Mummert, \emph{Reverse mathematics of {MF} spaces}, J. Math. Log. \textbf{6}
  (2006), no.~2, 203--232. \MR{2317427 (2008d:03011)}

\bibitem{MummertStephan10}
C.~Mummert and F.~Stephan, \emph{Topological aspects of poset spaces}, Michigan
  Math. J. \textbf{59} (2010), no.~1, 3--24. \MR{2654139}

\bibitem{Post44}
E.~L. Post, \emph{Recursively enumerable sets of positive integers and their
  decision problems}, Bull. Amer. Math. Soc. \textbf{50} (1944), 284--316.
  \MR{0010514 (6,29f)}

\bibitem{Simpson09}
S.~G. Simpson, \emph{Subsystems of second order arithmetic}, second ed.,
  Perspectives in Logic, Cambridge University Press, Cambridge, 2009.
  \MR{2517689 (2010e:03073)}

\bibitem{Smorynski77}
C.~Smory{\'n}ski, \emph{The incompleteness theorems}, Handbook of mathematical
  logic, Stud. Logic Foundations Math., vol.~90, North-Holland, Amsterdam,
  1977, pp.~821--865. \MR{0491063 (58 \#10343)}

\bibitem{Tychonoff35}
A.~Tychonoff, \emph{\"{U}ber die {A}bbildungen bikompakter {R}\"aume in
  {E}uklidische {R}\"aume}, Math. Ann. \textbf{111} (1935), no.~1, 760--761.
  \MR{1513029}

\bibitem{Watnick84}
R.~Watnick, \emph{A generalization of {T}ennenbaum's theorem on effectively
  finite recursive linear orderings}, J. Symbolic Logic \textbf{49} (1984),
  no.~2, 563--569. \MR{745385 (85i:03152)}

\end{thebibliography}

\end{document}